\documentclass[a4paper,11pt]{article}

\usepackage{bbm}
\usepackage{amscd}
\usepackage{stmaryrd}
\usepackage{amsfonts}
\usepackage{mathrsfs}
\usepackage{latexsym,amssymb,amsmath,amscd,amscd,amsthm,amsxtra,xypic}
\usepackage[dvips]{graphicx}
\usepackage[all]{xy}
\usepackage{ulem}

\usepackage[toc,page]{appendix}

\usepackage{geometry}
\usepackage{hyperref}

\usepackage{color}

\newcommand{\dlam}{_\lambda}

\newcommand{\g}{\mathfrak{g}}
\newcommand{\m}{\mathfrak{m}}
\newcommand{\h}{\mathfrak{h}}
\newcommand{\D}{\mathrm{D}}
\newcommand{\degree}{\mathrm{degree}}
\newcommand{\Der}{\mathrm{Der}}
\newcommand{\Courant}[1]{\left\llbracket  #1\right\rrbracket }

\newcommand{\G}{\mathcal{G}}

\newcommand{\C}{\mathcal{C}}

\DeclareMathOperator{\Hom}{Hom}

\DeclareMathOperator{\img}{img}
\DeclareMathOperator{\coker}{coker}

\DeclareMathOperator{\Diff}{Diff}

\DeclareMathOperator{\gl}{gl}
\DeclareMathOperator{\Aut}{Aut}

\DeclareMathOperator{\ad}{ad}

\DeclareMathOperator{\aff}{aff}
\DeclareMathOperator{\Mor}{Mor}

\newtheorem{Thm}{Theorem}[section]

\newtheorem{Pro}[Thm]{Proposition}
\newtheorem{Lem}[Thm]{Lemma}

\newtheorem{Def-Pro}[Thm]{Definition-Proposition}
\newtheorem{Def}[Thm]{Definition}

\theoremstyle{definition}
\newtheorem{Ex}[Thm]{Example}
\newtheorem{Rm}[Thm]{Remark}

\begin{document}
\title{A review of Lie $2$-algebras}

\author{\vspace{2mm}Honglei Lang  \,  and \, Zhangju Liu \\ \vspace{2mm}
\it{Department of Applied Mathematics, China Agricultural University, Beijing, 100083, China}\\ \vspace{2mm}
\it{Department of Mathematics, Peking University, Beijing, 100871, China}\\ \vspace{2mm}
hllang@cau.edu.cn,~~liuzj@pku.edu.cn}


\date{}
\maketitle

\makeatletter
\newif\if@borderstar
\def\bordermatrix{\@ifnextchar*{%
\@borderstartrue\@bordermatrix@i}{\@borderstarfalse\@bordermatrix@i*}%
}
\def\@bordermatrix@i*{\@ifnextchar[{\@bordermatrix@ii}{\@bordermatrix@ii[()]}}
\def\@bordermatrix@ii[#1]#2{%
\begingroup
\m@th\@tempdima8.75\p @\setbox\z@\vbox{%
\def\cr{\crcr\noalign{\kern 2\p@\global\let\cr\endline }}%
\ialign {$##$\hfil\kern 2\p@\kern\@tempdima & \thinspace %
\hfil $##$\hfil && \quad\hfil $##$\hfil\crcr\omit\strut %
\hfil\crcr\noalign{\kern -\baselineskip}#2\crcr\omit %
\strut\cr}}%
\setbox\tw@\vbox{\unvcopy\z@\global\setbox\@ne\lastbox}%
\setbox\tw@\hbox{\unhbox\@ne\unskip\global\setbox\@ne\lastbox}%
\setbox\tw@\hbox{%
$\kern\wd\@ne\kern -\@tempdima\left\@firstoftwo#1%
\if@borderstar\kern2pt\else\kern -\wd\@ne\fi%
\global\setbox\@ne\vbox{\box\@ne\if@borderstar\else\kern 2\p@\fi}%
\vcenter{\if@borderstar\else\kern -\ht\@ne\fi%
\unvbox\z@\kern-\if@borderstar2\fi\baselineskip}%
\if@borderstar\kern-2\@tempdima\kern2\p@\else\,\fi\right\@secondoftwo#1 $%
}\null \;\vbox{\kern\ht\@ne\box\tw@}%
\endgroup
}
\makeatother
\begin{abstract}
We first recall two equivalent definitions of  Lie $2$-algebras, categorification of Lie algebras and $2$-term $L_\infty$-algebras. Then we present four different kinds of Lie $2$-algebras from $2$-plectic manifolds, Courant algebroids, homotopy Poisson manifolds and affine multivector fields on a Lie groupoid respectively. Moreover, we recall the cohomology theory of Lie $2$-algebras and analyze its lower degree cases. The integration of strict Lie $2$-algebras to strict Lie $2$-groups is also discussed.

\end{abstract}

\tableofcontents

\section{Introduction}
As categorification of Lie algebras, the notion of Lie $2$-algebras was introduced by Baez and Crans in 2004 \cite{BC}. It is one of the fundamental objects in higher Lie theory and has close connection with strongly homotopy Lie algebras (also called $L_\infty$-algebras) introduced by  Schlessinger and Stasheff in \cite{SS}. See also  \cite{LM,LSt}. A Lie $2$-algebra is given by 
replacing the Jacobi identity in a Lie algebra with an isomorphism, called the Jacobiator, which satisfies a new law of its own.  Likewise, if the associative law in a Lie group is substituted with an isomorphism, called the associator, one gets a Lie $2$-group; see \cite{BL}.
Since its appearance, the Lie $2$-algebra structure has enjoyed significant applications in both geometry and mathematical physics. Its own algebraic properties have also drawn people's attention.

The purpose of this paper is to give a brief overview of Lie $2$-algebras on certain aspects. We first present the definition from two equivalent points of view, the categorification and the $2$-term $L_\infty$-algebra.  Explicitly, a Lie $2$-algebra is given by a bilinear skew-symmetric bracket on a linear category, which is a functor and obeys the Jacobi identity only up to a trilinear coherent natural transformation (the Jacobiator). It is equivalent to a $2$-term $L_\infty$-algebra, which is a $2$-term chain complex of vector spaces with a $2$-bracket and a $3$-bracket. The Jacobi identity of the $2$-bracket fails but is controlled by the $3$-bracket which satisfies its own law. In particular, when the Jacobiator is an identity, i.e., the $3$-bracket is trivial, one gets a strict Lie $2$-algebra. Strict Lie $2$-algebras are equivalent to Lie algebra crossed modules, which are classified by the third cohomology of a Lie algebra \cite{Gerstenhaber}. When the source and target of any morphism in  a Lie $2$-algebra are equal, i.e., the differential of the $2$-term complex is zero, we get a skeletal Lie $2$-algebra. Skeletal Lie $2$-algebras one-to-one correspond to triples that consist of a Lie algebra, a Lie algebra module and  a $3$-cocycle.  In fact, every Lie $2$-algebra is equivalent to a skeletal one. Thus  Lie $2$-algebras are classified by using the third cohomology of Lie algebras. See \cite{BC}.

Then we review several Lie $2$-algebras that come from geometric structures, namely, $2$-plectic manifolds, Courant algebroids, homotopy Poisson manifolds and affine multivector fields on Lie groupoids. A $2$-plectic structure on a manifold is a nondegenerate closed $3$-form. Similar to the Poisson bracket on the functions of a symplectic manifold, there is  a Lie $2$-algebra structure on functions and Hamiltonian $1$-forms of a $2$-plectic manifold \cite{BHR,R1,R2}. This Lie $2$-algebra is further used to define homotopy  moment maps \cite{CFRZ, FLZ, MZ}. A Courant algeboid \cite{LWX} is a vector bundle together with a bilinear form, a skew-symmetric bracket and an anchor map. The bracket satisfies the Jacobi identity up to a coboundary, which generates a Lie $2$-algebra on the spaces of sections of the bundle and functions on the base manifold \cite{RW}. Parallel to the fact  that there is a one-to-one correspondence between Lie algebra structures on a vector space and linear Poisson structures on the dual space, there is a one-to-one correspondence between Lie $2$-algebra structures on  a pair of vector spaces and linear homotopy Poisson structures on the dual \cite{LSX}.  
A homotopy Poisson structure on a graded manifold  is an $L_\infty$-algebra structure on the functions which is compatible with the graded commutative algebra structure by the Leibniz rule; see \cite{CF, KV, LSX, Vo}.  At last, on a Lie groupoid, multivector fields that are compatible with the groupoid multiplication are called multiplicative. Multiplicative multivector fields with the Schouten bracket form a graded Lie algebra, which is not invariant  under the Morita equivalence of Lie groupoids. Thus, to define multivector fields on a differentiable stack, one needs to extend the Lie algebra to a Lie $2$-algebra formed by affine multivector fields on a Lie groupoid, which is Morita invariant; see \cite{BLL,BCLX,LLS2,OW}.

The cohomology theory of a Lie $2$-algebra with a representation on a $2$-vector space is built up and lower degree cases are studied.  The first cocycles with respect to the adjoint representation on itself is used to define derivations of a Lie $2$-algebra \cite{CSZ, LLS}. A derivation Lie $2$-algebra and a derivation Lie $3$-algebra are further constructed in order to classify  nonabelian extensions of a Lie $2$-algebra \cite{CSZ}. The second cohomology is to study deformations of a Lie $2$-algebra and it classifies the abelian extension of a Lie $2$-algebra \cite{LSZ}. The third cohomology is to classify strong crossed modules for Lie $2$-algebras \cite{LL}. The corresponding results for Lie algebras are well-known.

At last, we review the integration of strict Lie $2$-algebras  to strict Lie $2$-groups, which is the same to integrating a Lie algebra crossed module to a Lie group crossed module. In particular, the integration of the derivation Lie $2$-algebra is the automorphism Lie $2$-group \cite{LLS}. 

The review is far from complete. Due to the limitations of our knowledge and the space of the paper, many important applications of Lie $2$-algebras  are missed here. For example, (1): as Lie algebras and Lie groups describe the symmetries of particles, Lie $2$-algebras and Lie $2$-groups are used to illustrate symmetries of symmetries. So Lie $2$-algebras play essential roles in higher gauge theory \cite{BS, BSCS, Waldorf}. In fact, connections on principal $2$-bundles are defined as Lie $2$-algebra valued $1$-forms; (2):  the integration of  a nonstrict Lie $2$-algebra is not discussed in detail here. We refer to \cite{BC, BSCS,Ge,H,SZ3,SZ1}; (3): Lie $2$-bialgebras  (\cite{BSZ, CSX2,Olga}) and the integration to Poisson Lie $2$-groups  are missed here; see \cite{CSX}.

This paper is organized as follows. Section $2$ is about the two equivalent definitions of Lie $2$-algebras and some basic facts about strict and skeletal Lie $2$-algebras. A Lie $2$-algebra from the skew-symmetrization of a Leibniz algebra is also given. In Section $3$, we review four Lie $2$-algebra structures from $2$-plectic manifolds, Courant algebroids, homotopy Poisson manifolds and affine multivector fields on a Lie groupoid, respectively. In section $4$, we recall the cohomology theory of Lie $2$-algebras and analyze the significance of lower degree cohomologies. The last section is about the integration of strict Lie $2$-algebras to strict Lie $2$-groups. We discuss the integration of the derivation Lie $2$-algebra to the automorphism Lie $2$-group.  Appendix A and B include a brief introduction on strict $2$-categories and Lie groupoids and algebroids with some examples.

\section{Definitions and examples}\label{definition and example}
In this section, we first introduce Lie $2$-algebras from the categorification point of view and by explicit formulas, under the names of 
semistrict Lie $2$-algebras \cite{BC} and $2$-term $L_\infty$-algebras \cite{SS} respectively.  See \cite{LM, LSt} for more about $L_\infty$-algebras. Then two special cases, strict Lie $2$-algebras and skeletal Lie $2$-algebras, are discussed.  Strict Lie $2$-algebras are equivalent to Lie algebra crossed modules, which can be classified by the $3$rd cohomology of Lie algebras \cite{Gerstenhaber}. 
Skeletal Lie $2$-algebras are used to classify Lie $2$-algebras  in terms of the $3$rd cohomology of Lie algebras. At last, a class of Lie $2$-algebras  is constructed by the skew-symmetrization of  Leibniz algebras.
\subsection{Lie $2$-algebras and $2$-term $L_\infty$-algebras}\label{2.1}
A Lie $2$-algebra is the categorification of a Lie algebra. It is a linear category with a functorial skew-symmetric bilinear operation satisfying the Jacobi identity up to natural transformations which obey coherence laws of their own.
To get a Lie $2$-algebra, we first categorify a vector space to a $2$-vector space. Let $\mathrm{Vect}$ denote the category of vector spaces. A {\bf $2$-vector space} is a category in $\mathrm{Vect}$. Explicitly,  a $2$-vector space is a category whose objects and morphisms form vector spaces and all the structure maps are linear.


\begin{Def} \rm{(\cite{BC})}
A {\bf (semistrict) Lie $2$-algebra} is a $2$-vector space $L$ equipped with
\begin{itemize}
\item  a skew-symmetric bilinear functor, the {\bf bracket}, $[\cdot,\cdot]:L\times L\to L$;
\item  a completely skew-symmetric trilinear natural isomorphism, the {\bf Jacobiator},
\[J_{x,y,z}:[[x,y],z]\to [x,[y,z]]+[[x,z],y],\]
that is required to satisfy the {\bf Jacobiator identity}:
\begin{eqnarray*}
&&J_{[w,x],y,z}([J_{w,x,z},y]+1)(J_{w,[x,z],y}+J_{[w,z],x,y}+J_{w,x,[y,z]})=\\ && [J_{w,x,y},z](J_{[w,y],x,z}+J_{w,[x,y],z})([J_{w,y,z},x]+1)([w,J_{x,y,z}]+1),
\end{eqnarray*}
\end{itemize}
for all $w,x,y,z\in L_0$, where $L_0$ is the object space of $L$. It is called  a {\bf strict} Lie $2$-algebra if the Jacobiator is the identity and a {\bf skeletal} Lie $2$-algebra if the source and target of any morphism are equal.
\end{Def}
\begin{Rm}The Jacobiator identity can be interpreted by a commutative octagon; see \cite{BC} for the diagram. 
\end{Rm}

The bracket in a semistrict Lie $2$-algebra is strictly skew-symmetric but obeys the Jocobi identity only up to a coherent natural transformation. If both the skew-symmetry and the Jacobi identity are allowed to hold up to natural transformations, we arrive at a notion called {\bf weak Lie $2$-algebras} introduced by Rotenberg in \cite{RLie2}. Weak Lie $2$-algebras are the complete categorification of Lie algebras.  In particular, weak Lie $2$-algabras with strict skew-symmetrizator are called {\bf semistrict} and with trivial Jacobiator are called {\bf hemistrict}. Semistrict and hemistrict Lie $2$-algebras form full sub-2-categories of weak Lie $2$-algebras. 
Throughout this paper,  we only study semistrict Lie $2$-algebras and call them simply as  Lie $2$-algebras.

Lie $2$-algebras are equivalent to an earlier notion called $2$-term $L_\infty$-algebras. $L_\infty$-algebras, also called strongly homotopy Lie algebras, were introduced by Stasheff and Schlessinger \cite{SS} in the study of deformation theory. 

\begin{Def}
A {\bf $2$-term $L_\infty$-algebra} $\mathfrak{g}$ is a $2$-term chain complex of vector spaces $\mathfrak{g}_{-1}\xrightarrow{d}\mathfrak{g}_0$ with 
\begin{itemize}
\item two skew-symmetric bilinear maps $l_2^0: \mathfrak{g}_0\wedge \mathfrak{g}_0\to \mathfrak{g}_0$ and $l_2^1: \mathfrak{g}_0\wedge\mathfrak{g}_{-1}\to \mathfrak{g}_{-1}$, which are always denoted as $l_2$ or $[\cdot,\cdot]$;
\item a totally skew-symmetric trilinear map $l_3: \mathfrak{g}_0\wedge \mathfrak{g}_0\wedge \mathfrak{g}_0\to \mathfrak{g}_{-1}$,
\end{itemize}
satisfying that
\begin{itemize}
\item[\rm{(1)}] $d[x,a]=[x,da],\quad [da,b]=[a,db]$;
\item[\rm{(2)}] $[[x,y],z]+[[y,z],x]+[[z,x],y]=-dl_3(x,y,z)$;
\item[\rm{(3)}] $[[x,y],a]+[[y,a],x]+[[a,x],y]=-l_3(x,y,da)$;
\item[\rm{(4)}] $l_3([w,x],y,z)-l_3([w,y],x,z)+l_3([w,z],x,y)+l_3([x,y],w,z)+l_3([y,z],w,x)-l_3([x,z],w,y)-
[l_3(w,x,y),z]-[l_3(w,y,z),x]+[l_3(w,x,z),y]+[l_3(x,y,z),w]=0$,
\end{itemize}
for all $w,x,y,z\in \mathfrak{g}_0$ and $a,b\in \mathfrak{g}_{-1}$. \end{Def}

To build the equivalence of Lie $2$-algebras and $2$-term $L_\infty$-algebras, we first relate $2$-vector spaces with
$2$-term chain complexes of vector spaces.
  A {\bf $2$-term chain complex} of vector spaces is a pair of vector spaces with a linear map between them, i.e.,  $C_{-1}\xrightarrow{d} C_0$.
Given a $2$-vector space $\mathbbm{V}$, denote by $V_{-1}$ and $V_0$ the spaces of morphisms and objects, and $s,t:V_{-1}\to V_0$ the source and target maps. Define $d:=t|_{\ker s}$, then 
\[\ker s \xrightarrow{d} V_0\]
is a $2$-term chain complex of vector spaces. Conversely, for a $2$-term chain complex $C_{-1}\xrightarrow{d} C_0$, define a $2$-vector space $\mathbbm{V}$ by
\[V_{-1}=C_{-1}\oplus C_0,\qquad V_0=C_0. \] The source, target maps and the composition of $\mathbbm{V}$ are 
\[s(\overrightarrow{f},x)=x,\qquad t(\overrightarrow{f},x)=x+d\overrightarrow{f},\qquad (\overrightarrow{f},x)\cdot (\overrightarrow{g},y)=(\overrightarrow{f}+\overrightarrow{g},y),\qquad \]
for $(\overrightarrow{f},x),(\overrightarrow{g},y)\in C_{-1}\oplus C_0$ with $x=y+d\overrightarrow{g}$. 

Actually, $2$-vector spaces, together with linear functors between them and linear natural transformations between functors, form a (strict) $2$-category, denoted by $\mathrm{2Vect}$. Similarly, by $\mathrm{2Term}$, we mean the $2$-category whose objects are $2$-term chain complexes of vector spaces, morphisms are chain maps and $2$-morphisms are chain homotopies between chain maps. 
By \cite[Proposition 8]{BC}, the $2$-category $\mathrm{2Vect}$ is $2$-equivalent to the $2$-category $\mathrm{2Term}$.
Denote by $\mathrm{Lie2Alg}$ and $\mathrm{2TermL_\infty}$ the $2$-categories of Lie $2$-algebras and $2$-term $L_\infty$-algebras respectively. 
We refer to Appendix A for details of $2$-categories.

\begin{Thm} \rm{(\cite{BC})}
The $2$-categories $\mathrm{Lie2Alg}$ and $\mathrm{2TermL_\infty}$ are $2$-equivalent.
\end{Thm}

As the description of $2$-term $L_\infty$-algebras is more explicit and convenient for calculation, throughout this paper, unless particular emphasis,  a Lie $2$-algebra is pointed to a $2$-term $L_\infty$-algebra.  
\begin{Def}\label{homo}
Let $\mathfrak{g}$ and $\mathfrak{h}$ be Lie $2$-algebras. A {\bf Lie $2$-algebra (or $2$-term $L_\infty$-algebra) homomorphism} $\phi:\mathfrak{g}\to \mathfrak{h}$ consists of 
\begin{itemize}
\item a chain map $(\phi_0,\phi_1):\mathfrak{g}\to \mathfrak{h}$, namely, $\phi_0: \mathfrak{g}_0\to \mathfrak{h}_0$ and $\phi_1:\mathfrak{g}_{-1}\to \mathfrak{h}_{-1}$ such that $\phi_0\circ d_\g=d_\h\circ \phi_1$,
\item a skew-symmetric bilinear map $\phi_2: \mathfrak{g}_0\wedge \mathfrak{g}_0\to \mathfrak{h}_{-1}$,
such that the following equations hold for $x,y,z\in \mathfrak{g}_0$ and $a\in \mathfrak{g}_{-1}$:
\end{itemize}
\begin{itemize}
\item[\rm{(1)}] $\phi_0[x,y]_\g-[\phi_0(x),\phi_0(y)]_\h=d_\h\phi_2(x,y)$,
\item[\rm{(2)}] $\phi_1[x,a]_\g-[\phi_0(x),\phi_1(a)]_\h=\phi_2(x,d_\g a)$,
\item[\rm{(3)}] $\phi_1l_3^\g(x,y,z)-l_3^\h(\phi_0(x),\phi_0(y),\phi_0(z))=[\phi_0(x),\phi_2(y,z)]_\h+[\phi_0(y),\phi_2(z,x)]_\h+[\phi_0(z),\phi_2(x,y)]_\h-\phi_2([x,y]_\g,z)-\phi_2([y,z]_\g,x)-\phi_2([z,x]_\g,y)$.
\end{itemize}
If $\phi_2=0$, it is called a {\bf strict} Lie $2$-algebra homomorphism.
\end{Def}

Clearly, $\mathrm{Id}=(\mathrm{id}_{\g_0},\mathrm{id}_{\g_{-1}},0):\g\longrightarrow\g$ is a Lie $2$-algebra homomorphism, which is called the {\bf identity homomorphism}. Let $\phi:\g\to \g'$ and $\psi:\g'\to \g''$ be two Lie 2-algebra homomorphisms. Their composition $\psi\diamond \phi:\g\to \g''$ is also a homomorphism, which is defined by $(\psi\diamond \phi)_0=\psi_0\circ \phi_0:\g_0\to \g''_0$, $(\psi\diamond \phi)_1=\psi_1\circ \phi_1:\g_{-1}\to \g''_{-1}$
and
$$(\psi\diamond \phi)_2=\psi_2\circ (\phi_0\times \phi_0)+\psi_1\circ \phi_2:\g_0\wedge \g_0\to \g''_{-1},$$
A Lie $2$-algebra homomorphism $\phi:\g\to \g'$ is called an {\bf isomorphism} if there exists a Lie $2$-algebra homomorphism $\phi^{-1}:\g'\to \g$ such that
the compositions $\phi^{-1}\diamond \phi:\g\to \g$ and $\phi\diamond  \phi^{-1}:\g'\to \g'$ are both the identity homomorphism. In particular,  an {\bf automorphism} of a Lie $2$-algebra $\g$  is a Lie $2$-algebra isomorphism from $\g$ to itself. 

\begin{Lem}
  Let $\phi=(\phi_0,\phi_1,\phi_2):\g\to \g'$ be
a Lie $2$-algebra homomorphism. If $\phi_0$ and $\phi_1$ are invertible, then $\phi$ is an isomorphism, and  $\phi^{-1}$  is given by
 $$\phi^{-1}=(\phi^{-1}_0, \phi^{-1}_1, -\phi^{-1}_1 \circ \phi_2\circ (\phi^{-1}_0\times \phi^{-1}_0)).$$
\end{Lem}

\subsection{Strict and skeletal Lie $2$-algebras}
A Lie $2$-algebra $(\g,d,[\cdot,\cdot],l_3)$ is called {\bf strict} if $l_3=0$. Then the $2$-bracket $[\cdot,\cdot]$  satisfies the Jacobi identity. For a $2$-vector space (or, a $2$-term chain complex), its functors and natural transformations of a $2$-vector space constitute a strict Lie $2$-algebra.
\begin{Ex}Let $\mathbbm{V}: V_{-1}\xrightarrow{\partial} V_0$ be a $2$-term chain complex of vector spaces. Its truncated morphisms constitute a strict Lie $2$-algebra
\[\gl(\mathbbm{V}): \gl_{-1}(\mathbbm{V}):=\Hom(V_0,V_{-1})\xrightarrow{\delta} \gl_0(\mathbbm{V}),\qquad \delta(D)=(\partial\circ D,D\circ \partial),\]
where
\[\gl_0(\mathbbm{V}):=\{(A_0,A_1)\in \gl(V_0)\oplus \gl(V_{-1});A_0\circ \partial=\partial\circ A_1\},\]
and the $2$-bracket is defined by the commutator:
\begin{eqnarray*}
{}[(A_0,A_1),(B_0,B_1)]_C&=&(A_0\circ B_0-B_0\circ A_0,A_1\circ B_1-B_1\circ A_1),\\  {}[(A_0,A_1),D]_C&=&A_1\circ D-D\circ A_0,
\end{eqnarray*}
for $(A_0,A_1),(B_0,B_1)\in \gl_0(\mathbbm{V})$ and $D\in \Hom(V_0,V_{-1})$. We see that $\gl_0(\mathbbm{V})$ is the space of functors from $\mathbbm{V}$ to $\mathbbm{V}$ and $\gl_{-1}(\mathbbm{V})$ is the space of natural transformations.

\end{Ex}
This strict Lie $2$-algebra will be used to define representations of Lie $2$-algebras; see Subsection \ref{com}.

An equivalent description of strict Lie $2$-algebras is the concept of  Lie algebra crossed modules, which first appeared in the work of Gerstenhaber \cite{Gerstenhaber}. 

\begin{Def}\label{Liecm}
A {\bf Lie algebra crossed module} is a quadruple $(\m,\g,\varphi,\triangleright)$ consisting of Lie algebras $\m$ and $\g$, a Lie algebra homomorphism $\varphi:\m\to \g$, and an action $\triangleright$ of $\g$ as derivations of $\m$ (that is, a Lie algebra homomorphism $\triangleright:\g \to \Der(\m)$), satisfying 
\begin{eqnarray*}
\varphi(x\triangleright a)&=&[x,\varphi(a)];\\ 
\varphi(a)\triangleright b&=&[a,b],
\end{eqnarray*}
for $x\in \g$ and $a,b\in \m$.
\end{Def}

A Lie algebra crossed module becomes a strict Lie $2$-algebra if taking the action of $\g$ on $\m$ as the bracket. Conversely, given a strict Lie $2$-algebra $\g_{-1}\xrightarrow{d} \g_0$, it is obvious that $\g_{-1}$ with the bracket $[a,b]=[da,b]$ is a Lie algebra and the Lie algebra $\g_0$ acts on $\g_{-1}$ by $x\triangleright a=[x,a]$. Then $\g$ is a Lie algebra crossed module.  By this analysis, we have
\begin{Pro}
There is a one-to-one correspondence between strict Lie $2$-algebras and Lie algebra crossed modules.
\end{Pro}

\begin{Ex}
Let $\g$ be a Lie algebra and $\h\subset \g$ an ideal. Then $\h\hookrightarrow \g$ is a Lie algebra crossed module, where the action of $\g$ on $\h$ is given by the Lie bracket on $\g$.

\end{Ex}
\begin{Ex}
Let $\g$ be a Lie algebra and $\Der(\g)$  its derivation Lie algebra. Then $\ad:\g\to \Der(\g)$ defined by $\ad(x)=[x,\cdot]$ is a Lie algebra crossed module, where $\Der(\g)$ acts on $\g$ naturally.
\end{Ex}

It was first shown by Gerstenhaber \cite{Gerstenhaber} that Lie algebra crossed modules are classified by the third cohomology of a Lie algebra. A direct proof can be found in \cite{Wa}. 
Given a Lie algebra crossed module $\m\xrightarrow{\varphi} \g$, there exists a $4$-term exact sequence
\begin{equation*}
\xymatrix@C=0.5cm{ 0 \ar[r] & V\ar[rr]^{i} &&
                \mathfrak{m} \ar[rr]^{\varphi} && \mathfrak{g} \ar[rr]^{\pi} && \mathfrak{h}\ar[r]  & 0,
                }
\end{equation*}
where the cokernel $\mathfrak{h}$ is a Lie algebra and the kernel
$V$ is an $\mathfrak{h}$-module induced by the action of
$\g$ on $\m$.

\begin{Def}
Two crossed modules $\varphi:\m\to \g$ and $\varphi':\m'\to \g'$ such that $\ker \varphi=\ker \varphi'=V$ and $\coker \varphi=\coker \varphi'=\h$ are called {\bf elementary equivalent} if there are Lie algebra homomorphisms $\phi_1:\m\to \m'$ and $\phi_0:\g\to \g'$ such that $\phi_1(x\triangleright a)=\phi_0(x)\triangleright' \phi_1(a)$ and the following diagram is commutative:
\begin{equation*}
\CD
  0 @> >>  V @> >> \mathfrak{m} @>\varphi>> \mathfrak{g}@> >> \mathfrak{h} @> >> 0 \\
 @. @V \mathrm{id} VV @V \phi_1VV @V \phi_0 VV @V \mathrm{id} VV @.  \\
 0 @> >> V @> >> \mathfrak{m'} @>\varphi'>> \mathfrak{g'}@>  >> \mathfrak{h} @> >> 0.
\endCD
\end{equation*}
\end{Def}

Define the {\bf equivalence relation} on Lie algebra crossed modules by the equivalence relation generated by elementary equivalences.
Denote by crmod$(\mathfrak{h},V)$ the set of
equivalence classes of crossed modules with fixed kernel
$V$, fixed cokernel $\mathfrak{h}$ and fixed action of $\h$ on $V$. 
\begin{Thm} \rm{(Gerstenhaber)}\label{Gerstenhaber} With notations above,
there is a bijection
\[ \mathrm{crmod} (\mathfrak{h},V) \cong 
\mathrm{H}^3(\mathfrak{h},V).\]
\end{Thm}

Another special case  is the skeletal Lie $2$-algeba.
A Lie $2$-algebra $(\g,d,[\cdot,\cdot],l_3)$ is called {\bf skeletal} if $d=0$. Then the $2$-bracket $[\cdot,\cdot]$ satisfies the Jacobi identity. Thus $(\g_0,[\cdot,\cdot])$ is a Lie algebra, $\g_{-1}$ becomes a $\g_0$-module, and $l_3:\wedge^3 \g_0\to \g_{-1}$ is a $3$-cocycle.

\begin{Ex}\label{string}
Given a semisimple Lie algebra $\mathfrak{g}$ with the Killing form $(\cdot,\cdot)$, 
we have a skeletal Lie $2$-algebra
\[\mathbbm{R}\xrightarrow{0} \mathfrak{g},\]
where the $2$-bracket $[\cdot,\cdot]$ is the Lie bracket on $\mathfrak{g}$ and trivial otherwise, and $l_3(x,y,z)=(x,[y,z])$.  Note that $l_3$ is a $3$-cocycle of $\g$ with coefficients in $\mathbbm{R}$, i.e., $l_3\in Z^3(\mathfrak{g},\mathbbm{R})$. This Lie $2$-algebra is called the {\bf string Lie $2$-algebra} and denoted by $\mathrm{String}(\g)$.
\end{Ex}


\begin{Pro}\rm{(\cite{BC})}\label{ske}
There is a one-to-one correspondence between isomorphism classes of skeletal Lie $2$-algebras and isomorphism classes of quadruples consisting of a Lie algebra $\g$, a vector space $V$, a representation of $\g$ on $V$, and an element in the third cohomology $\mathrm{H}^3(\g,V)$.\end{Pro}

For a Lie $2$-algebra $\g_{-1}\xrightarrow{d} \g_0$, write $\g_0$ and $\g_{-1}$ as: $\g_0=\img d\oplus \g'_0$ and $\g_{-1}=\ker d\oplus \g'_{-1}$. This allows us to define a $2$-term chain complex with trivial differential:
\[\ker d\xrightarrow{0} \g'_0.\]
Actually, as shown in \cite{BC}, we have 
\begin{Lem}
Every Lie $2$-algebra is equivalent, as an object of $\mathrm{Lie2Alg}$, to a skeletal one.
\end{Lem}
This result together with Proposition \ref{ske} leads to a classification of all Lie $2$-algebras up to equivalence.

\begin{Thm}\rm{(\cite{BC})}\label{classification}
There is a one-to-one correspondence between equivalence classes of Lie $2$-algebras 
and isomorphism classes of quadruples as in Proposition \ref{ske}.
\end{Thm}

For the string Lie $2$-algebra in Example \ref{string}, since $\g$ is semi-simple, the Killing form is non-degenerate. So the $3$-cocycle $l_3$ described above represents a nontrivial cohomology class. 
By Theorem \ref{classification}, the Lie $2$-algebra $\mathrm{String(\g)}$ 
 is not equivalent to a skeletal strict Lie $2$-algebra.


\subsection{Leibniz algebras and Lie $2$-algebras}\label{lebniz}
The skew-symmetrization of a Leibniz algebra gives a Lie $2$-algebra  \cite{LS}. In fact, a Leibniz algebra gives rise to a hemistrict Lie $2$-algebra and the skew-symmetrization of any hemistrict Lie $2$-algebra is a semistrict Lie $2$-algebra; see \cite{RLie2}.
Besides this, the skew-symmetrization of a Leibniz algebra can be extended in a natural way to a Lie-Yamaguti algebra structure; see \cite{KW}.  

A {\bf Leibniz algebra} is a vector space $\g$ endowed with a bilinear bracket $\{\cdot,\cdot\}:\g\times \g\to \g$ satisfying the Leibniz rule:
\[\{x,\{y,z\}\}=\{\{x,y\},z\}+\{y,\{x,z\}\},\qquad \forall x,y,z\in \g.\]
The {\bf left center} of a Leibniz algebra $\g$ is defined by
\[Z(\g)=\{x\in \g;\{x,y\}=0,\forall y\in \g\}.\]
$Z(\g)$ is obviously an ideal of $\g$.
Denote by $\Courant{\cdot,\cdot}$  the skew-symmetrization of the Leibniz bracket
\begin{eqnarray*}\label{Cbracket}
\Courant{x,y}=\frac{1}{2}(\{x,y\}-\{y,x\}).
\end{eqnarray*}
By $J_{x,y,z}$, we mean the corresponding Jacobiator, i.e., 
\[J_{x,y,z}=\Courant{x,\Courant{y,z}}+\Courant{y, \Courant{z,x}}+\Courant{z, \Courant{x,y}},\qquad \forall x,y,z\in \g.\]
It is shown that  $J_{x,y,z}\in Z(\g)$.

\begin{Thm}\rm(\cite{LS})\label{Leibniz}
Let $(\g,\{\cdot,\cdot\})$ be a Leibniz algebra and $Z(\g)$ its left center. Then we have a Lie $2$-algebra
\[Z(\g)\hookrightarrow \g,\]
where 
\begin{equation*}
\left\{\begin{array}{rll}
l_2(x,y)&=&\Courant{x,y}=\frac{1}{2}(\{x,y\}-\{y,x\});\\
l_2(x,c)&=&\Courant{x,c}=\frac{1}{2}\{x,c\};\\
l_3(x,y,z)&=&J_{x,y,z},
\end{array}\right.
\end{equation*}
for all $x,y,z\in \g$ and $c\in Z(\g)$.
\end{Thm}
An omni-Lie algebra is a typical example of Leibniz algebras, which was introduced in \cite{W} by linearizing Courant algebroids. Its Dirac structures characterize all the Lie algebra structures on a vector space.
An {\bf omni-Lie algebra} is a triple $(\gl(V)\oplus V, \langle\cdot,\cdot\rangle, \{\cdot,\cdot\})$, where $V$ is a vector space, $\langle\cdot,\cdot\rangle$ is a nondegenerate pairing 
\[\langle A+u,B+v\rangle=Av+Bu,\qquad \forall A+u,B+v\in \gl(V)\oplus V,\]
and the bracket is given by
\[\{A+u,B+v\}=[A,B]_C+Av,\]
where $[\cdot,\cdot]_C$ is the commutator bracket. 
Omni-Lie algebras are Leibniz algebras, i.e.,
\[\{e_1,\{e_2,e_3\}\}=\{\{e_1,e_2\},e_3\}+\{e_2,\{e_1,e_3\}\},\qquad \forall e_1,e_2,e_3\in \gl(V)\oplus V.\]
By Theorem \ref{Leibniz}, we have the following example of Lie $2$-algebras from an omni-Lie algebra.
\begin{Ex}\label{omni}
For a vector space $V$,  there is a natural Lie $2$-algebra 
\[V\hookrightarrow \gl(V)\oplus V,\]
where 
\begin{equation*}
\left\{\begin{array}{rll}
l_2^0(A+u,B+v)&=&[A,B]+\frac{1}{2}(Av-Bu);\\
l_2^1(A+u,v)&=&\frac{1}{2}Av;\\
l_3(A+u,B+v,C+w)&=&-\frac{1}{4}([A,B]w+[B,C]u+[C,A]v),
\end{array}\right.
\end{equation*}
for all $A+u,B+v,C+w\in \gl(V)\oplus V$. Here $[\cdot,\cdot]$ is the commutator bracket $[\cdot,\cdot]_C$.
\end{Ex}

\section{Lie $2$-algebras from geometric structures}
Lie $2$-algebras appear in various geometric structures. In this part, we shall review four geometric structures, $2$-plectic manifolds, Courant algebroids, homotopy Poisson manifolds and affine multivector fields on Lie groupoids that give rise to Lie $2$-algebras. Some applications are also given. We refer to Appendix B for a brief introduction of Lie groupoids and algebroids.

\subsection{Lie $2$-algebras from $2$-plectic manifolds}  

In the classical mechanics of point particles, the phase space is often a symplectic manifold and the poisson bracket of its functions makes the space of  observables a Lie algebra. Passing from point particles to strings, the phase space becomes a $2$-plectic manifold and a Lie $2$-algebra of observables is constructed. This Lie $2$-algebra is used to describe the dynamics of a classical bosonic string. See \cite{BHR}.

\begin{Def}
A  $3$-form $\omega\in \Omega^3(M)$ on a smooth manifold $M$ is called a {\bf $2$-plectic structure} if it is both closed:
\[d_{dR}\omega=0,\]
and nondegenerate, i.e., for $X\in T_m M$, 
\[ \omega(X,\cdot,\cdot)=0\Rightarrow  X=0.\]
We call the pair $(M,\omega)$ a {\bf $2$-plectic manifold}.
\end{Def}

In classical mechanics, the position and momentum of a particle are given by a point in $T^*M$, which is a symplectic manifold with the symplectic structure given by the differential of a canonical $1$-form. To describe the position and momentum of a string, we need $\wedge^2 T^*M$, which is a $2$-plectic manifold with the $3$-form given by the differential of a canonical $2$-form.

\begin{Ex}
Let $M$ be a manifold. On $\wedge^2 T^*M$, there is a canonical $2$-form $\alpha\in \Omega^2(\wedge^2 T^*M)$ given by:
\[\alpha_\xi(v_1,v_2)=\xi(\pi_*(v_1),\pi_*(v_2)),\qquad \forall \xi\in \wedge^2 T^*M,v_1,v_2\in T_\xi (\wedge^2 T^*M),\]
where $\pi: \wedge^2 T^*M\to M$ is the projection.  Then $\wedge^2 T^*M$ with the $3$-form 
\[\omega=d_{dR} \alpha\in \Omega^3(\wedge^2 T^*M)\]
is $2$-plectic manifold.
\end{Ex}

Let $(M,\omega)$ be a $2$-plectic manifold. From the nondegeneracy of $\omega$, we have an injective map
\[T_m M\to \wedge^2 T_m^* M,\qquad X\mapsto \omega(X,\cdot,\cdot),\]
which is not an isomorphism in general.
\begin{Def}
Let $(M,\omega)$ be a $2$-plectic manifold. A $1$-form $\alpha$ on $M$ is called {\bf Hamiltonian} if there exists a vector field $X_\alpha$
on $M$ such that 
\[d_{dR}\alpha=\omega(X_\alpha,\cdot,\cdot).\]
$X_\alpha$ is called the {\bf Hamiltonian vector field} corresponding to $\alpha$.
\end{Def}
Denote by $\Omega^1_{\mathrm{Ham}}(M)$ the vector space of Hamiltonian $1$-forms. 
Similar to the Poisson bracket on functions in the symplectic case, there is a bracket on Hamiltonian $1$-forms:
\[\{\alpha,\beta\}=\omega(X_\alpha,X_\beta,\cdot),\qquad \forall \alpha,\beta\in  \Omega^1_{\mathrm{Ham}}(M).\]

\begin{Thm}\rm{(\cite{BHR,R1})}\label{2-plectic}
If $(M,\omega)$ is a $2$-plectic manifold, there is a  Lie $2$-algebra 
\[L_\infty(M,\omega): C^\infty(M)\xrightarrow{d_{dR}}\Omega^1_{\mathrm{Ham}}(M),\]
where $d_{dR}$ is the de Rham differential, and the brackets are 
\begin{equation*}
\left\{\begin{array}{rll}
l_2(\alpha,\beta)&=&\{\alpha,\beta\};\\
l_2(\alpha,f)&=&0;\\
l_3(\alpha,\beta,\gamma)&=&\omega(X_\alpha,X_\beta,X_\gamma),
\end{array}\right.
\end{equation*}
for $\alpha,\beta,\gamma\in \Omega^1_{\mathrm{Ham}}(M)$ and $f\in C^\infty(M)$.
\end{Thm}
This result can be generalized to an arbitrary $n$. An $n$-plectic manifold is a manifold with a nondegenerate closed $(n+1)$-form, which is canonically equipped with a Lie $n$-algebra structure \cite{R1}.

The following example is related to the Wess-Zumino-Witten model and loop groups.
\begin{Ex}
Let $G$ be a compact simple Lie group with Lie algebra $\mathfrak{g}$ and let $(\cdot,\cdot)$ be the Killing form on $\mathfrak{g}$. There is a $2$-plectic structure $\omega\in \Omega^3(G)$ on $G$  such that
\[\omega(x,y,z)=(x,[y,z]),\qquad \forall x,y,z\in T_e G=\g,\]
and  $\omega$  is both left and right-invariant.

Let $\Omega^1_{\mathrm{Ham}}(G)^L$ be the set of left-invariant Hamiltonian $1$-forms. Actually, $\Omega^1_{\mathrm{Ham}}(G)^L=\mathfrak{g}^*$.
Since the left-invariant smooth functions on $G$ are constants, we have a Lie $2$-algebra $L_\infty(G,\omega)$ on the complex:
\[0: \mathbbm{R}\to \mathfrak{g}^*.\]
This Lie $2$-algebra is in fact the string Lie $2$-algebra $\mathrm{String}(\g)$ in Example \ref{string}.
\end{Ex}

The Lie $2$-algebra $L_\infty(M,\omega)$ associated with a $2$-plectic manifold $(M,\omega)$  is used to define Lie algebra and Lie $2$-algebra moment maps \cite{CFRZ, FLZ, MZ}. 
See a Lie algebra $\g$  as a strict Lie $2$-algebra $0\to \g$.
\begin{Def} \rm{(\cite{CFRZ})} 
Let $\g\to \mathfrak{X}(M),u\mapsto \hat{u}$ be a Lie algebra action on a $2$-plectic manifold $(M,\omega)$ such that $\hat{u}$ for all $u\in \g$ are Hamiltonian vector fields. A {\bf homotopy moment map} for this action (or a {\bf $\g$-moment map}) is a Lie $2$-algebra homomorphism
\[(\phi_0,\phi_1=0,\phi_2):\g\to L_\infty(M,\omega),\]
where $\phi_0:\g\to \Omega^1_{\mathrm{Ham}}(M)$ and  $\phi_2:\wedge^2 \g\to C^\infty(M)$, such that
\[-\iota_{\hat{u}} \omega=d_{dR}(\phi_0(u)),\qquad \forall u\in \g.\]
\end{Def}

Fixing a point $p\in M$, define $\omega_{3p}:\wedge^3 \g\to \mathbbm{R}$ by 
\[\omega_{3p}(u_1,u_2,u_3):=\omega(\hat{u_1},\hat{u_2},\hat{u_3})|_p,\qquad \forall u_1,u_2,u_3\in \g.\]
It can be checked that $\omega_{3p}$ is a $3$-cocycle of the Lie algebra $\g$ with coefficients in $\mathbbm{R}$, i.e., $\omega_{3p}\in Z^3(\g,\mathbbm{R})$. Moreover, the cohomology class $[\omega_{3p}]\in \mathrm{H}^3(\g,\mathbbm{R})$ is independent of the choice of the point $p$.

\begin{Thm}\rm{(\cite{CFRZ})}
The existence of a $\g$-moment map implies that $[\omega_{3p}]=0$. The converse is true if $\mathrm{H}^1(M)=0$.

\end{Thm}
In \cite{MZ}, the authors went one step further by replacing the Lie algebra $\g$ by a skeletal Lie $2$-algebra. We would not go into the details here.

\subsection{Lie $2$-algebras from Courant algebroids}

Courant algebroids were originally introduced in \cite{LWX} for studying the double of Lie bialgebroids.  They include as examples the double of Lie bialgebras and the bundle $TM\oplus T^*M$ with the bracket introduced by Courant \cite{C} for the study of Dirac structures. Courant algebroids play important roles in generalized complex geometry \cite{G}. See also \cite{D,R, SW} for other applications.
Courant algebroids give rise to Lie $2$-algebras naturally \cite{RW}.

\begin{Def} \rm{(\cite{LWX})}
A {\bf Courant algebroid} is a vector bundle $E\to M$ equipped with 
\begin{itemize}
\item a nondegenerate symmetric bilinear form $\langle\cdot,\cdot\rangle$ on the bundle;
\item  a skew-symmetric bracket $\Courant{\cdot,\cdot}$, called the Courant bracket, on $\Gamma(E)$;
\item  a bundle map $\rho:E\to TM$, called the anchor,
\end{itemize}
such that the following properties hold:
\begin{itemize}
\item[\rm{(1)}] $\Courant{\Courant{e_1,e_2},e_3}+\Courant{\Courant{e_2,e_3},e_1}+\Courant{\Courant{e_3,e_1},e_2}=\mathcal{D}T(e_1,e_2,e_3)$;
\item[\rm{(2)}] $\rho\Courant{e_1,e_2}=[\rho(e_1),\rho(e_2)]$;
\item[\rm{(3)}] $\Courant{e_1,fe_2}=f\Courant{e_1,e_2}+\rho(e_1)(f)e_2-\frac{1}{2}\langle e_1,e_2\rangle \mathcal{D} f$;
\item[\rm{(4)}]  $\rho\circ \mathcal{D}=0$, i.e., $\langle \mathcal{D}f,\mathcal{D}g\rangle=0$;
\item[\rm{(5)}] $\rho(e_1)\langle e_2,e_3\rangle=\langle \Courant{e_1, e_2}+\frac{1}{2}\mathcal{D}\langle e_1,e_2\rangle,e_3\rangle+\langle e_2,\Courant{e_1, e_3}+\frac{1}{2}\mathcal{D}\langle e_1,e_3\rangle\rangle$,
\end{itemize}
for $e_1,e_2,e_3\in \Gamma(E)$ and $f\in C^\infty(M)$. Here 
 $\mathcal{D}:C^\infty(M)\to \Gamma(E)$ is defined by 
 \begin{eqnarray}\label{DCA}
\langle \mathcal{D}f,e\rangle=\rho(e)(f),\qquad \forall e\in \Gamma(E), f\in C^\infty(M),
\end{eqnarray}
and $T(e_1,e_2,e_3)\in C^\infty(M)$ is defined by
\begin{eqnarray}\label{TCA}
T(e_1,e_2,e_3)=\frac{1}{6}(\langle \Courant{e_1,e_2},e_3\rangle+\langle \Courant{e_2, e_3},e_1\rangle+\langle \Courant{e_3, e_1},e_2\rangle).
\end{eqnarray}
\end{Def}
In \cite{LWX}, the authors also introduced a bracket:
\[\{e_1,e_2\}:=\Courant{e_1, e_2}+\mathcal{D}\langle e_1,e_2\rangle,\qquad \forall e_1,e_2\in \Gamma(E).\]
It is not skew-symmetric. 
Later, this bracket is found to satisfy the Jacobi identity
\[\{e_1, \{e_2,e_3\}\}=\{\{e_1,e_2\},e_3\}+\{e_2,\{e_1,e_3\}\},\]
and is called the {\bf Dorfman bracket}.
Rotenberg \cite{R} gives a simplified definition of Courant algebroids by using the Dorfman bracket. 


\begin{Thm}\rm{(\cite{RW})}\label{CALie2}
Let $(E,\langle \cdot,\cdot\rangle,\Courant{\cdot,\cdot},\rho)$ be a Courant algebroid. Then we obtain a Lie $2$-algebra
\[L_\infty(E): C^\infty(M)\xrightarrow{\mathcal{D}} \Gamma(E),\]
where $l_2$ and $l_3$ are given by
\begin{equation*}
\left\{\begin{array}{rll}
l_2(e_1,e_2)&=&\Courant{e_1,e_2};\\
l_2(e_1,f)&=&\frac{1}{2}\langle e_1,\mathcal{D}f\rangle;\\
l_3(e_1,e_2,e_3)&=&-T(e_1,e_2,e_3),
\end{array}\right.
\end{equation*}
for $e_1,e_2,e_3\in \Gamma(E)$ and $f\in C^\infty(M)$, where $\mathcal{D}$ and $T$ are defined in \eqref{DCA} and \eqref{TCA} respectively.
\end{Thm}
More generally, from a pre-Courant algebroid, one is also able to get a Lie $2$-algebra; see \cite{LiuSX}.

\begin{Ex}\label{Standard}
On the bundle $TM\oplus T^*M$, define a bilinear form by the pairing 
\[\langle X+\alpha,Y+\beta\rangle=\alpha(Y)+\beta(X),\qquad \forall X,Y\in \mathfrak{X}(M),\alpha,\beta\in \Omega^1(M),\]
 the Courant bracket by
\begin{eqnarray}\label{CTT}
\Courant{X+\alpha,Y+\beta}=[X,Y]+L_X \beta-L_Y \alpha+\frac{1}{2}d_{dR}(\alpha(Y)-\beta(X)),
\end{eqnarray}
and the anchor $\rho$ by the projection to $TM$.
Then $(TM\oplus T^*M,\langle \cdot,\cdot\rangle,\Courant{\cdot,\cdot},\rho)$ is a Courant algeboid. It was used in \cite{C} to study Dirac structures.
\end{Ex}

In particular, let $V$ be a vector space and $M=V^*$. Choosing linear vector fields and constant forms on $V^*$, one gets $\gl(V)\oplus V$, which inherits the pairing, Courant bracket and the anchor from the Courant algebroid $TV^*\oplus T^*V^*$. It turns out to be the omni-Lie algebra we recalled in Subsection \ref{lebniz}. So omni-Lie algebras are the linearization of Courant algebroids. In fact, the Lie $2$-algebra in Example \ref{omni} is a particular case of Theorem \ref{CALie2}.

It was then developed in \cite{SW} that the Courant bracket \eqref{CTT} on $TM\oplus T^*M$ can be twisted by a $3$-form $\omega\in \Omega^3(M)$ in the way:
\[\Courant{X+\alpha,Y+\beta}_\omega=\Courant{X+\alpha,Y+\beta}-\iota_Y\iota_X \omega.\]
\begin{Ex}
Let $M$ be a manifold with a $3$-form $\omega\in \Omega^3(M)$. Then $TM\oplus T^*M$ with the bilinear form and anchor the same as in Example \ref{Standard} and the bracket $\Courant{\cdot,\cdot}_\omega$ is a Courant algebroid if and only if $d_{dR}\omega=0$.
\end{Ex}

By Theorem \ref{CALie2}, from a manifold $M$ with a closed $3$-form $\omega\in \Omega^3(M)$, we get a Lie $2$-algebra 
$C^\infty(M)\xrightarrow{d_{dR}} \mathfrak{X}(M)\oplus \Omega^1(M)$. Denote it by $L_\infty(TM\oplus T^*M)$.
Recall from Theorem \ref{2-plectic} that we have another Lie $2$-algebra $L_\infty(M,\omega)$ if $\omega\in \Omega^3(M)$ is closed and nondegenerate. The relation of these two Lie $2$-algebras is clarified in \cite{R2} as follows.

\begin{Thm}
With the notations above, there exists a Lie $2$-algebra homomorphism $(\phi_0,\phi_1,\phi_2)$ embedding $L_\infty(M,\omega)$ into $L_\infty(TM\oplus T^*M)$:
\begin{equation*}
\xymatrix{ C^\infty(M) \ar@{->} [r]^{\phi_1=\mathrm{id}}\ar@ {->} [d]_{d_{dR}}
&C^\infty(M) \ar@ {->} [d]^{d_{dR}}\\
{\Omega^1_{\mathrm{Ham}}(M)} \ar@{->} [r]^{\phi_0 }&\mathfrak{X}(M)\oplus \Omega^1(M),}
\end{equation*}
where $\phi_0(\alpha)=(X_\alpha,\alpha)$ with $X_\alpha$ being the Hamiltonian vector field of $\alpha$, i.e., $d_{dR}\alpha=\iota_{X_\alpha} \omega$, and 
$\phi_2:\wedge^2\Omega^1_{\mathrm{Ham}}(M)\to C^\infty(M)$ is given by
\[\phi_2(\alpha,\beta)=\frac{1}{2}(\alpha(X_\beta)-\beta(X_{\alpha})),\qquad \forall \alpha,\beta\in \Omega^1_{\mathrm{Ham}}(M).\]
\end{Thm}

We refer to \cite{R2} for the discussion of this result for exact Courant algebroids.

\subsection{Lie $2$-algebras from homotopy Poisson manifolds}
For a vector space $V$, there is a one-to-one correspondence between Lie algebra structures on $V$ and linear Poisson structures on $V^*$. Such a result for Lie $2$-algebras was explored in \cite{LSX}. In fact, on the dual of a Lie $2$-algebra, 
there is a linear homotopy Poisson structure and also a linear quasi-Poisson groupoid structure. Moreover, from a Lie $2$-algebra, there associates a natural Courant algebroid structure.

A homotopy Poisson algebra is a graded commutative algebra with an $L_\infty$-structure
 whose brackets satisfy the Leibniz rule. It has appeared as higher Poisson structures in \cite{KV, Vo} and
$P_\infty$-structures in \cite{CF}.

\begin{Def}
A {\bf  homotopy Poisson algebra} of degree $n$   is a graded commutative algebra $\mathfrak{a}$ with an $L_\infty$-algebra structure $\{l_m\}_{m\geq1}$ on $\mathfrak{a}[n]$, such that the map
\[x\longrightarrow{l_m(x_1,\cdots,x_{m-1},x)},\qquad \forall x_1,\cdots,x_{m-1},x\in \mathfrak{a}\]
is a derivation of degree $2-m-n(m-1)+\sum_{i=1}^{m-1}|x_i|$.
Here, $|x|$ denotes the degree of $x\in\mathfrak{a}.$

A  homotopy Poisson algebra of degree $n$ is of {\bf finite type} if there exists a $q$ such that $l_m=0$ for all $m>q$.
A {\bf homotopy Poisson manifold} of degree $n$ is a graded manifold $\mathcal{M}$ whose algebra of functions $C^\infty(\mathcal{M})$ is equipped with a degree $n$ homotopy Poisson algebra structure of finite type.
\end{Def}

Given a Lie $2$-algebra $(\g,d, [\cdot,\cdot],l_3)$, its graded dual space $\g^*[1]=\mathfrak{g}^*_0[1]\oplus{\mathfrak{g}^*_ {-1}[1]}$ is an N-manifold of degree $1$ with the base manifold $\g_{-1}^*$.  Its algebra of  functions is
$$
\cdots\oplus \big(C^\infty(\g_{-1}^*)\otimes \wedge^2\g_{0}\big)\oplus \big(C^\infty(\g_{-1}^*)\otimes\g_{0}\big)\oplus C^\infty(\g_{-1}^*).
$$
There is a degree $1$ homotopy Poisson algebra structure on it obtained by extending the original Lie $2$-algebra structure using the Leibniz rule. Thus, the dual of a Lie $2$-algebra is a linear homotopy Poisson manifold of degree $1$. Here linear means that the brackets of linear functions are still linear functions.
The converse is easy to see.

\begin{Pro}
Let $\mathbbm{V}=V_{-1}\oplus V_0$ be a graded vector space. Then there is a one-to-one correspondence between Lie $2$-algebra structures on $\mathbbm{V}$ and linear homotopy Poisson structures on $\mathbbm{V}^*[1]$.
\end{Pro}

Let $(\g,d,[\cdot,\cdot],l_3)$ be a Lie $2$-algebra. 
The cotangent bundle $T^*[2]\g^*[1]$ is a symplectic NQ-manifold of degree $2$. Symplectic NQ-manifolds of degree $2$ are in one-to-one
correspondence with Courant algebroids \cite{R3}. Thus we obtain a Courant algebroid $E$ from a Lie $2$-algebra $\g$:
\begin{equation*}\label{eq:E}
 E=
\g_{-1}^*\times (\g_0^*\oplus \g_0)\longrightarrow
 \g_{-1}^*,
\end{equation*}
in which the anchor and the Dorfman bracket are defined by the derived bracket using the degree $3$ function $l=-d-l_2-l_3$. Here elements in $\g_0,\g_0^*,\g_{-1},\g_{-1}^*$ as functions are of degree  $1,1,0,2$ respectively.

\begin{Pro}\rm{(\cite{LSX})} The Courant algebroid structure on $E$ is as follows: the  bilinear form is the pairing between $\g_0$ and $\g_0^*$; the anchor and Dorfman bracket are determined by  \begin{itemize}
  \item[\rm(i)]$\rho(x)(a)=[x,a]$, the anchor of $x$ is a linear vector field on $\g_{-1}^*$;
      \item[\rm(ii)]$\rho(\xi)=d^*(\xi)$, the anchor of $\xi$ is a constant vector field on $\g_{-1}^*$;
    \item[\rm(iii)]  $\{x,y\}=[x,y]+l_3(x,y,\cdot)\in\mathfrak{g}_{0}+\g_0^*\otimes\mathfrak{g}_{-1};$
    \item[\rm(iv)] $\{\xi,\eta\}=0$ and $\{x,\xi\}=-\ad_x^*\xi$, where $\ad_x^*:\g_0^*\to \g_0^*$ is given by $\langle \ad_x^*\xi,y\rangle=-\langle \xi,[x,y]\rangle$,
            \end{itemize}
  for constant sections $x,y\in{\mathfrak{g}_0}, \xi,\eta\in{\mathfrak{g}_{0}^*}$ and linear function $a\in{\mathfrak{g}_{-1}}$.
          
          \end{Pro}
   So a Courant algebroid is derived from a Lie $2$-algebra, which further gives rise to a new Lie $2$-algebra as in Theorem \ref{CALie2}. See \cite{LSX} for the relation of these two Lie $2$-algebras.   
 This Courant algebroid $E$ is actually the double of a Lie quasi-bialgebroid; see \cite{LWX}. By integration, we shall further get a quasi-Poisson groupoid from a Lie $2$-algebra.

On a Lie groupoid $\mathcal{G}$, a bivector field $\Pi$ is called {\bf multiplicative} if the graph of the groupoid multiplication 
$\Lambda= \{(g,h,gh);s(g)=t(h)\}\subset \mathcal{G}\times \mathcal{G}\times \mathcal{G}$ is coisotropic relative to $\Pi\oplus \Pi\oplus -\Pi$.  Namely, 
\[(\Pi\times \Pi\times -\Pi)(\xi_1,\xi_2)=0,\qquad \forall \xi_1,\xi_2\in \Lambda^\perp,\]
where 
\[\Lambda^\perp=\{\xi\in T^*_\lambda (\mathcal{G}\times \mathcal{G}\times \mathcal{G});\lambda\in \Lambda, \langle\xi, v\rangle=0, \forall v\in T_\lambda \Lambda\}.\]

\begin{Def} \rm{(\cite{Xu})}\label{quasi-Poisson}
A {\bf quasi-Poisson groupoid} is a Lie groupoid $\mathcal{G}$ with a multiplicative bivector field $\Pi\in \mathfrak{X}^2(\mathcal{G})$ and a $3$-section $\phi\in \Gamma(\wedge^3 A)$, where $A$ is the Lie algebroid of $\mathcal{G}$, such that
\[\frac{1}{2}[\Pi,\Pi]=\overleftarrow{\phi}-\overrightarrow{\phi},\qquad [\Pi,\overleftarrow{\phi}]=0.\]
If $\phi=0$, it is called a {\bf Poisson groupoid}.
\end{Def}

For the  $2$-term chain complex $\g_0^*\xrightarrow{d^*} \g_{-1}^*$, the abelian Lie group $\g_0^*$ acts on $\g_{-1}^*$ by 
\[\xi \triangleright \alpha=\alpha+d^*(\xi),\qquad \forall \xi\in \g_0^*,\alpha\in \g_{-1}^*.\]
This gives rise to an action Lie groupoid  $\mathcal{G}: \mathfrak{g}^*_{0}\times \mathfrak{g}^*_{-1}\rightrightarrows{\mathfrak{g}^*_{-1}}$, where the source, target and multiplication are
\[s(\xi,\alpha)=\alpha,\qquad t(\xi,\alpha)=\alpha+d^*(\xi),\qquad (\xi_1,\alpha_1)(\xi_2,\alpha_2)=(\xi_1+\xi_2,\alpha_2),\]
where $\alpha_1=\alpha_2+d^*(\xi_2)$.

The Lie algebroid of this Lie groupoid is the action Lie algebroid $A: \mathfrak{g}^*_{0}\times \mathfrak{g}^*_{-1}\to{\mathfrak{g}^*_{-1}}$. See Example \ref{actiongroupoid} in Appendix B for action Lie groupoids and algebroids. Observe that an element $l_3\in \wedge^3\g_0^*\otimes \g_{-1}$ can be viewed as a $3$-section of $A$.

\begin{Thm} \rm{(\cite{LSX})} Let $\g$ be a Lie $2$-algebra. Then we get a linear quasi-Poisson groupoid $(\mathcal{G}, \Pi,\phi)$, 
 where $\mathcal{G}: \mathfrak{g}^*_{0}\times \mathfrak{g}^*_{-1}\rightrightarrows{\mathfrak{g}^*_{-1}}$ is the action Lie groupoid described above, $\phi=-l_3\in \Gamma(\wedge^3 A)$ and $\Pi$ is characterized by
{\footnotesize
\[\Pi(d_{dR}x,d_{dR}y)=-[x,y],\qquad \Pi(d_{dR}x,d_{dR}a)=-[x,a], \qquad \Pi(d_{dR}a,d_{dR}b)=-[da,b],\]
}
where $x,y\in \g_0$ and $a,b\in \g_{-1}$ are  linear functions on $\mathfrak{g}^*_{0}\times{\mathfrak{g}^*_{-1}}$.
\end{Thm}
In summary, we have 
$$
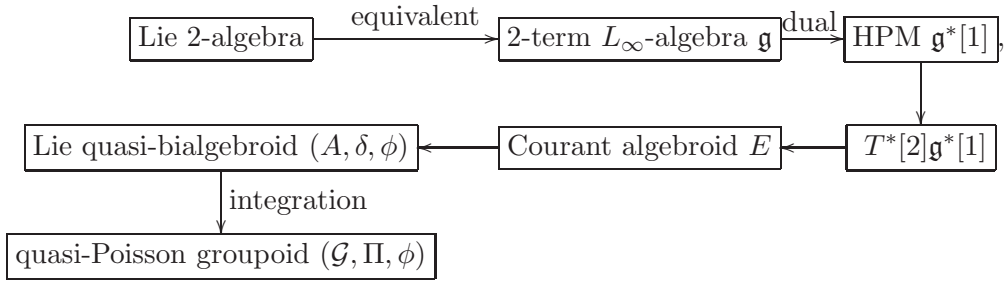
\xymatrix{
*+[F]{\mbox{Lie~ 2-algebra}}\ar[r]^{\mbox{equivalent}\quad} &*+[F]{\mbox{2-term~} L_\infty\mbox{-algebra}~\g}\ar[r]^{\mbox{\qquad dual}}&*+[F]{\mbox{HPM}~\g^*[1]}\ar[d] \\
*+[F]{\mbox{Lie~quasi-bialgebroid}~(A,\delta,\phi)}\ar[d]^{\mbox{integration}}&*+[F]{\mbox{Courant~algebroid}~E}\ar[l]&*+[F]{\mbox{}~T^*[2]\g^*[1]}\ar[l]\\
*+[F]{\mbox{quasi-Poisson~groupoid}~(\mathcal{G},\Pi,\phi)}.&&
 },
$$
where HPM stands for homotopy Poisson manifolds.
\subsection{Lie $2$-algebras from affine structures}
Geometric structures on a Lie groupoid that are compatible with the groupoid multiplication are called  multiplicative structures. For example,
a $k$-vector field $\Pi$ on a Lie groupoid $\mathcal{G}$ is called {\bf multiplicative} (\cite{Xu}) if the graph of the groupoid multiplication $\{(g,h,gh);s(g)=t(h)\}\subset \mathcal{G}\times \mathcal{G}\times \mathcal{G}$ is coisotropic with respect to $\Pi\oplus \Pi\oplus (-1)^{k+1}\Pi$. This generalizes the multiplicative bivector field for a quasi-Poisson groupoid in Definition \ref{quasi-Poisson}.

In \cite{BLL,OW},  the authors constructed a strict Lie $2$-algebra on the multiplicative $1$-vector fields on a Lie groupoid and their natural transformations.  This construction is Morita invariant and gives rise to a strict Lie $2$-algebra structure on the differentiable stack. Later this idea is generalized to multiplicative multivector fields in \cite{BCLX}. A strict graded Lie $2$-algebra is obtained which is used to define  multivector fields on a differentiable stack. A geometric explanation of this Lie $2$-algebra is given in \cite{LLS2} by using affine multivector fields. See \cite{Luthesis, Waffine} for affine structures on a Lie group.
In this subsection, we always suppose that $\mathcal{G}$ is a Lie groupoid with Lie algebroid $A$.
\begin{Def}\label{aff vf}
A $k$-vector field $\Pi\in \mathfrak{X}^k(\mathcal{G})$ on a Lie groupoid $\mathcal{G}$ is called {\bf affine} if the submanifold \[S:=\{(g,h,l,hg^{-1}l); s(g)=s(h),t(g)=t(l)\}\subset \mathcal{G}\times \mathcal{G}\times \mathcal{G}\times \mathcal{G}\]  is coisotropic with respect to $\Pi\oplus (-1)^{k+1}\Pi\oplus (-1)^{k+1} \Pi\oplus \Pi$.
\end{Def}
It is shown in
 \cite{Xu} that a $k$-vector field $\Pi\in \mathfrak{X}^k(\mathcal{G})$ is multiplicative if and only if it is affine and the base manifold $M$ is coisotropic with respect to $\Pi$. 
 
 \begin{Lem}\rm{(\cite{CSX})}\label{affine vf}
A $k$-vector field $\Pi\in \mathfrak{X}^k(\mathcal{G})$ on a Lie groupoid $\mathcal{G}$ is affine if and only if the following two conditions hold:
\begin{itemize}
\item[\rm(i)]
for any $(g,h)\in \mathcal{G}^{(2)}$,
\begin{eqnarray*}\label{affine}
\Pi(gh)=L_{\mathcal{X}*} \Pi(h)+R_{\mathcal{Y}*} \Pi(g)-L_{\mathcal{X}*}\circ R_{\mathcal{Y}*}(\Pi(s(g))),
\end{eqnarray*}
where $\mathcal{X}$ and $\mathcal{Y}$ are any two local bisections passing through $g$ and $h$ respectively.
\item[\rm(ii)] for any $\xi\in \Omega^1(M)$, $\iota_{t^*\xi} \Pi$ is right-invariant.
\end{itemize}
\end{Lem}

 \begin{Lem}\label{leftright} Let $\Pi$ be a  $k$-vector  field on the Lie groupoid $\mathcal{G}$ with $\pi=\mathrm{pr}_{\wedge^k A} \Pi|_M\in \Gamma(\wedge^k A)$.
Define
\begin{eqnarray}\label{affmul}
\Pi_r=\Pi-\overrightarrow{\pi},\qquad \Pi_l=\Pi-\overleftarrow{\pi}.
\end{eqnarray}
Then $\Pi$ is affine if and only if $\Pi_l$  or $\Pi_r$ is a multiplicative $k$-vector field on $\mathcal{G}$.
\end{Lem}

Denote by  $\mathfrak{X}^k_{\mathrm{aff}}(\mathcal{G})$ and $\mathfrak{X}^k_{\mathrm{mult}}(\mathcal{G})$ the spaces of affine and multiplicative $k$-vector fields respectively. We have $\mathfrak{X}^k_{\mathrm{mult}}(\mathcal{G})\subset \mathfrak{X}^k_{\mathrm{aff}}(\mathcal{G})$. The following result is from \cite{LLS2}.

\begin{Thm}\begin{itemize}
\item[\rm{(i)}]
We have a $2$-vector space $\mathfrak{X}^k_{\mathrm{aff}}(\mathcal{G})$, whose spaces of morphisms and objects are 
$\mathfrak{X}^k_{\mathrm{aff}}(\mathcal{G})$ and $\mathfrak{X}^k_{\mathrm{mult}}(\mathcal{G})$ respectively, 
and the source and target maps are given by $s(\Pi)=\Pi_r$ and $t(\Pi)=\Pi_l$ as defined in \eqref{affmul}. 
\item[\rm{(ii)}] The graded $2$-vector space
$\oplus_k \mathfrak{X}^k_{\aff}(\mathcal{G})$ with the Schouten bracket is a strict graded Lie $2$-algebra.
\end{itemize}
\end{Thm}
Here by a {\bf strict graded Lie $2$-algebra}, we say that the spaces of objects and morphisms are graded vector spaces and the Lie bracket is a graded Lie bracket. See \cite{BCLX} for the explicit definition.
The $2$-term $L_\infty$-algebra associated to this Lie $2$-algebra is the one in \cite{BCLX} stated as follows. Define $\Gamma(\wedge^\bullet A):=\oplus_k \Gamma(\wedge^k A)$ and $\mathfrak{X}^\bullet_{\mathrm{mult}}(\mathcal{G}):=\oplus_k\mathfrak{X}^k_{\mathrm{mult}}(\mathcal{G})$.

\begin{Pro}
There is a strict graded Lie $2$-algebra structure on 
\[ \Gamma(\wedge^\bullet A)\xrightarrow{d}  \mathfrak{X}^\bullet_{\mathrm{mult}}(\mathcal{G}),\] 
where $d(u)=\overrightarrow{u}-\overleftarrow{u}$, the bracket on $\mathfrak{X}^\bullet_{\mathrm{mult}}(\mathcal{G})$ is the Schouten bracket and the bracket $[\Pi,u]$ for  $\Pi\in \mathfrak{X}^k_{\mathrm{mult}}(\mathcal{G})$ and $u\in \Gamma(\wedge^l A)$ is determined by the relation $\overrightarrow{[\Pi,u]}=[\Pi,\overrightarrow{u}]$.
\end{Pro}
The homotopy equivalence class of this strict graded Lie $2$-algebra is invariant under the Morita equivalence of Lie groupoids, thus is considered as the space of multivector fields on the corresponding differentiable stack. We refer to Appendix A for the definition of homotopy equivalences of Lie $2$-algebras.

Let $\mathcal{G}\rightrightarrows M$ be a Lie groupoid and $J: X\to M$  a surjective smooth map. A {\bf left action} of $\mathcal{G}$ on $X$ along $J$, which is called the moment map, is a smooth map 
\[\mathcal{G}\times_M X=\{(g,x)\in \mathcal{G}\times X; s(g)=J(x)\}\to X,\qquad (g,x)\mapsto g\cdot x=gx\]
such that 
\[J(gx)=t(g),\qquad (gh)x=g(hx),\qquad 1_{J(x)} x=x.\]
We then say that  $X$ is a {\bf left $\mathcal{G}$-space}.
A {\bf left $\mathcal{G}$-bundle} is a left $\mathcal{G}$-space $X$ together with a $\mathcal{G}$-invariant surjective submersion $\beta: X\to N$. A left $\mathcal{G}$-bundle  is called {\bf principal} if the map
\[\mathcal{G}\times_M X\to X\times_N X,\qquad (g,x)\mapsto (gx,x)\]
is a diffeomorphism.

\begin{Def}
A {\bf Morita equivalence} between two Lie groupoids $\mathcal{G}_1\rightrightarrows M_1$ and $\mathcal{G}_2\rightrightarrows M_2$ is given by a principal $\mathcal{G}_1$-$\mathcal{G}_2$-bibundle, i.e., a manifold $X$ with moment maps $\alpha:X\to M_1$ and $\beta:X\to M_2$, such that $\beta: X\to M_2$ is a left principal $\mathcal{G}_1$-bundle, $\alpha:X\to M_1$ is a right principal $\mathcal{G}_2$-bundle and the two actions commute: $g_1\cdot (x\cdot g_2)=(g_1\cdot x)\cdot g_2$ for any $g_1\in \mathcal{G}_1,x\in X$ and $g_2\in \mathcal{G}_2$:
\begin{equation*}
\xymatrix{
              \mathcal{G}_1\ar@<0.5ex>[d] \ar@<-0.5ex>[d] \ar@/^1pc/ [r]&X  \ar[dl]_{\alpha} \ar[dr]^{\beta} & \mathcal{G}_2\ar@/_1pc/[l]\ar@<0.5ex>[d] \ar@<-0.5ex>[d]\\
M_1  &  &M_2}.
\end{equation*}
\end{Def}


\begin{Thm}\rm{(\cite{BCLX})}\label{BCLXth}
 Let $\mathcal{G}_1\rightrightarrows M_1$ and $\mathcal{G}_2\rightrightarrows M_2$  be Morita equivalent Lie groupoids. Then any 
 $\mathcal{G}_1$-$\mathcal{G}_2$-bibundle $M_1\leftarrow X\rightarrow M_2$ induces  a homotopy equivalence between the strict graded Lie $2$-algebras $\Gamma(\wedge^\bullet A_1)\to \mathfrak{X}^\bullet_{\mathrm{mult}}(\mathcal{G}_1)$ and 
$ \Gamma(\wedge^\bullet A_2)\to \mathfrak{X}^\bullet_{\mathrm{mult}}(\mathcal{G}_2)$.
\end{Thm}
\begin{Def}
 Let $\mathfrak{X}$ be a differentiable stack. The space of {\bf multivector fields} on $\mathfrak{X}$ is
defined to be the homotopy equivalence class of  strict graded Lie $2$-algebras $\Gamma(\wedge^\bullet A)\to \mathfrak{X}^\bullet_{\mathrm{mult}}(\mathcal{G})$, where $\mathcal{G}\rightrightarrows M$ is any Lie groupoid representing $\mathfrak{X}$.
 \end{Def}

In other words, a multivector field on a differentiable stack $\mathfrak{X}$ is defined as a homotopy equivalence class of affine multivector fields on any groupoid representing $\mathfrak{X}$.
Following \cite{BCLX}, a shifted $(+1)$ Poisson structure on a differentiable stack $\mathfrak{X}$ is simply  an element of the Maurer-Cartan moduli set of the  dgla decided by the homotopy equivalence class of strict graded Lie $2$-algebras corresponding to $\mathfrak{X}$.

\section{Modules and cohomology}
The cohomology theories for $L_\infty$-algebras and $A_\infty$-algebras were developed in \cite{LM,Ma,P} and also in \cite{Mi} for a bigger framework. In particular, the Lie $2$-algebra cohomology was formulated in \cite{BSZ} for the strict case to characterize strict Lie $2$-bialgebras and in \cite{LSZ} for the general to study deformations of Lie $2$-algebras. See also \cite{CLX} for the cohomology of hemistrict Lie $2$-algebras and \cite{A} for a cohomological theory of both strict Lie $2$-algebras and Lie $2$-groups. The $1$-cocycles of a Lie $2$-algebra with respect to the adjoint representation are used to define derivations for a Lie $2$-algebra. We shall also see that the second cohomology is used to classify the deformations and abelian extensions of Lie $2$-algebras \cite{LSZ}  and  the third cohomology is for classifying crossed modules of Lie $2$-algebras \cite{LL}. See \cite{D, Gerstenhaber, Wa} for the corresponding results for Lie algebras. 

\subsection{Definition}\label{com}
\begin{Def}
Let $\g$ be a Lie $2$-algebra and $\mathbbm{V}$ a $2$-vector space. A {\bf representation} of $\mathfrak{g}$ on  $\mathbbm{V}$ is a Lie $2$-algebra homomorphism $\phi:\mathfrak{g}\to \gl(\mathbbm{V})$. Such a $2$-vector space $\mathbbm{V}$ is called a {\bf $\mathfrak{g}$-module}.
\end{Def}
For simplicity, we always denote an action $\phi=(\phi_0,\phi_1,\phi_2):\mathfrak{g}\to \gl(\mathbbm{V})$ by
 $$x\triangleright (u+m):=\phi_0(x)(u+m),\qquad a\triangleright u:=\phi_1(a)(u),\qquad  (x,y)\triangleright u:=\phi_2(x,y)(u),$$
  for all  $x,y\in \g_0, a\in \g_{-1}, u\in V_0$ and $m\in V_{-1}.$

A Lie algebra $\g$ is a Lie $2$-algebra $0\to \g$. A {\bf $2$-term representation up to homotopy} of a Lie algebra $\g$ on a $2$-vector space $\mathbbm{V}$ is defined to be a Lie $2$-algebra homomorphism $\phi:\g\to \gl(\mathbbm{V})$ (\cite{SZ2}).
 \begin{Ex}\label{semidirect}
 Given a Lie algebra $(\g,[\cdot,\cdot]_\g)$ with a $2$-term representation up to homotopy on $\mathbbm{V}:V_{-1}\xrightarrow{\partial} V_0$, we get a Lie $2$-algebra
\[\g\ltimes \mathbbm{V}: V_{-1}\xrightarrow{0+\partial} \g\oplus V_0,\]
where the brackets are given by
\begin{equation*}
\left\{\begin{array}{rcl} 
l_2(x+u,y+v)&=&[x,y]_\g+x\triangleright v-y\triangleright u;\\ 
l_2(x+u,m)&=&x\triangleright m;\\ 
l_3(x+u,y+v,z+w)&=&-(x,y)\triangleright w-(y,z)\triangleright u-(z,x)\triangleright v,
\end{array}\right.
\end{equation*}
for $x,y,z\in \g, u,v,w\in V_0$ and $m\in V_{-1}$.
This Lie $2$-algebra is called the {\bf semidirect product Lie $2$-algebra}. We refer to \cite{SZ2} for more explicit examples.
\end{Ex}
\begin{Ex}\label{adjoint}
Every Lie $2$-algebra $(\g,d,[\cdot,\cdot],l_3)$ has a natural representation on itself, called the {\bf adjoint representation}, which is given by   $\ad=(\ad_0,\ad_1,\ad_2): \g\to \gl(\g)$, where  
\[\ad_0(x)=[x,\cdot],\qquad 
\ad_1(a)=[a,\cdot],\qquad \ad_2(x,y)=-l_3(x,y,\cdot),\]
for $x,y,z\in \g_0$ and $a\in \g_{-1}$.

\end{Ex}

 Let $\g$ be a Lie $2$-algebra and $\mathbb{V}:
V_{-1}\stackrel{\partial}{\longrightarrow}V_{0}$  a $\g$-module given by $\phi:\g\to \gl(\mathbbm{V})$.  
The cohomology of $\g$ comes from
the generalized Chevalley-Eilenberg complex as follows:
{\footnotesize
\begin{equation*} \label{CE complex}
\begin{split}
 &\degree \ -1: V_{-1}\stackrel{\D}{\longrightarrow}\\
 &\degree \ \ \ 0\ \ : V_{0}\oplus \Hom(\g_{0},V_{-1})\stackrel{\D}{\longrightarrow}\\
  &\degree\ \ \ 1\ \ : \Hom(\g_{0},V_{0})\oplus\Hom(\g_{-1},V_{-1})\oplus\Hom(\wedge^2\g_0,V_{-1})\stackrel{\D}{\longrightarrow}\\
  &\degree\ \ \ 2\ \ : \Hom(\g_{-1},V_{0})\oplus \Hom(\wedge^2\g_0,V_{0})\oplus \Hom(\g_0\wedge \g_{-1},V_{-1}) \oplus \Hom(\wedge^3\g_0,V_{-1})\stackrel{\D}{\longrightarrow}\\
  & \stackrel{\D}{\longrightarrow}\cdots,
\end{split}
\end{equation*}}
where the degrees of $\g_0,\g_{-1},V_0,V_{-1}$ are $-1,-2,0,-1$ respectively. Denote by $C^k(\g,\mathbb{V})$ the set of $k$-cochains. The coboundary
operator $\D$ can be decomposed as:
\begin{equation}\label{DD}
  \D=\hat{d}+\hat{\partial}+d^{(1,0)}_\phi+d^{(0,1)}_\phi+d_{\phi_2}+d_{l_3},
\end{equation}
where, for $s=0,-1$,
\begin{eqnarray*}
\hat{d}:\Hom(\wedge^p\g_0\wedge \odot
^q\g_{-1},V_s)&\longrightarrow&\Hom(\wedge  ^{p-1}\g_0\wedge
\odot  ^{q+1}\g_{-1},V_s),\\
\hat{\partial}:\Hom(\wedge^ p\g_0\wedge \odot
^q\g_{-1},V_{-1})&\longrightarrow&\Hom(\wedge^p\g_0\wedge
\odot ^q\g_{-1},V_{0}),\\
d_\phi^{(1,0)}:\Hom(\wedge  ^p\g_0\wedge \odot
^q\g_{-1},V_s)&\longrightarrow& \Hom(\wedge  ^{p+1}\g_0\wedge
\odot  ^q\g_{-1},V_s),\\
d_\phi^{(0,1)}:\Hom(\wedge  ^p\g_0\wedge \odot  ^q\g_{-1},V_0)&\longrightarrow&
\Hom(\wedge  ^p\g_0\wedge \odot  ^{q+1}\g_{-1},V_{-1}),\\
d_{\phi_2}:\Hom(\wedge^{p}\g_0\wedge \odot^q\g_{-1},V_0)&\longrightarrow& \Hom(\wedge^{p+2}\g_0\wedge \odot^q\g_{-1},V_{-1}),\\
d_{l_3}:\Hom(\wedge^p\g_0\wedge \odot^q \g_{-1},V_s)&\longrightarrow& \Hom(\wedge^{p+3}\g_0\wedge \odot^{q-1} \g_{-1},V_s).
\end{eqnarray*}
More concretely, for all
${x_i}\in{\mathfrak{g}_0},a_i\in{\mathfrak{g}_{-1}}, i\in{\mathbb{N}},$
{\footnotesize
\begin{eqnarray*}
&&\hat{d}f(x_1,\cdots,x_{p-1},a_1,\cdots,a_{q+1})
 =(-1)^{p}\big(f(x_1,\cdots,x_{p-1},da_1,a_2,\cdots,a_{q+1})+c.p.(a_1,\cdots,a_{q+1})\big),\\
&&\hat{\partial}f=(-1)^{p+2q}\partial\circ f,\\
&&d_\phi^{(1,0)}f(x_1,\cdots,x_{p+1},a_1,\cdots,a_{q})
 =\sum_{i=1}^{p+1}(-1)^{i+1}x_i \triangleright f(x_1,\cdots,\widehat{x_i},\cdots,x_{p+1},
  a_1,\cdots,a_{q})\\
  &&\ \ \ \ \ \ \ \ \ \ \ \ \ \ \ \ \ \ \ \ \ \ \ \ \ \ \ \ \ \ \ \ \ \ \ \ \ \ \ \ \ \ \ \ +\sum_{i<j}(-1)^{i+j}f([x_i,x_j],x_1,\cdots,\widehat{x_i},\cdots,\widehat{x_j},\cdots,x_{p+1},
  a_1,\cdots,a_{q})\\
  &&\ \ \ \ \ \ \ \ \ \ \ \ \ \ \ \ \ \ \ \ \ \ \ \ \ \ \ \ \ \ \ \ \ \ \ \ \ \ \ \ \ \ \ \ +\sum_{i,j}(-1)^{i}f(x_1,\cdots,\widehat{x_i},\cdots,x_{p+1},
  a_1,\cdots,[x_i,a_j],\cdots,a_{q}),\\
&&d_\phi^{(0,1)}f(x_1,\cdots,x_{p},a_1,\cdots,a_{q+1})
 =\sum_{i=1}^{q+1}(-1)^{p}a_i \triangleright f(x_1,\cdots,x_{p},
  a_1,\cdots,\widehat{a_i},\cdots,a_{q+1}),\\
&&d_{\phi_2}f(x_1,\cdots,x_{p+2},a_1,\cdots,a_{q})=
  \sum_{\sigma}(-1)^{p+2q}(-1)^\sigma(x_{\sigma(1)},x_{\sigma(2)})\triangleright f(x_{\sigma(3)},\cdots,x_{\sigma(p+2)},a_1,\cdots,a_{q}),\\
\label{lem:proerty p}&&d_{l_3}f(x_1,\cdots,x_{p+3},a_1,\cdots,a_{q-1})=\sum_\tau-(-1)^\tau f(x_{\tau(4)},\cdots,x_{\tau(p+3)},a_1,\cdots,a_{q-1},l_3(x_{\tau(1)},x_{\tau(2)},x_{\tau(3)})),
\end{eqnarray*}}
where $\sigma$ and $\tau$ are taken over all $(2,p)$-unshuffles and $(3,p)$-unshuffles respectively.
Define \[Z^k(\g,\mathbbm{V}):=\{f\in C^k(\g,\mathbbm{V});\D f=0\},\qquad B^k(\g,\mathbbm{V}):=\{\D g\in C^k(\g,\mathbbm{V}); g\in C^{k-1}(\g,\mathbbm{V})\}.\]
Elements in $Z^k(\g, \mathbbm{V})$ and $B^k(\g,\mathbbm{V})$ are called {\bf $k$-cocycles} and {\bf $k$-coboundaries}
respectively. The space $\mathrm{H}^k(\g,\mathbbm{V}):=Z^k(\g,\mathbbm{V})/B^k(\g,\mathbbm{V})$ is called the {\bf $k$th  cohomology} of the Lie $2$-algebra $\g$ with respect to the $\g$-module $\mathbbm{V}$.

\subsection{The first cohomology and derivation Lie $2$-algebras}
Let us write down the expressions of 1-cocycles and 1-coboundaries for a Lie $2$-algebra $\g$ with a $\g$-module $\mathbbm{V}$.
\begin{Lem}\label{1chain}
Let $(X,l_X )\in C^1(\g,\mathbb{V})$, where $X=(X_0,X_1) \in
\Hom(\g_0,V_0)\oplus\Hom(\g_{-1},V_{-1})$ and $l_X\in\Hom(\wedge^2
\g_0,V_{-1})$. Then
 \begin{itemize}
\item[\rm{(i)}]
 $(X,l_X)\in Z^1(\g, V)$ if and only if
\begin{eqnarray*}
\label{d} \partial\circ X_1&=&X_0\circ d,\\
\label{a}
\partial l_{X}(x,y)&=&X_0[x,y]+y\triangleright X_0x-x\triangleright X_0y,\\
\label{b}l_{X}(x,da)&=&X_1[x,a]+a\triangleright X_0x-x\triangleright X_1a,\\
\label{c}X_1l_{3}(x,y,z)&=&(l_{X}(x,[y,z] )+x\triangleright
l_{X}(y,z)-(y,z)\triangleright X_0x)+c.p.,
\end{eqnarray*}
where $x,y,z\in \g_0$ and $a\in \g_{-1}$.
\item[\rm{(ii)}] $(X,l_X)\in B^1(\g, V)$ if and only if $\exists (u,\Theta)\in V_0\oplus
\Hom(\g_0,V_{-1})=C^0(\g,\mathbb{V})$, s.t., 
\begin{eqnarray*}
\label{1-cob1}X(x+a)&=&x\triangleright u+a\triangleright u-\partial\Theta(x)-\Theta(da),\\
\label{1-cob2}l_X(x,y)&=&(x,y)\triangleright u+x\triangleright
\Theta(y)-y\triangleright \Theta(x)-\Theta([x,y]).
\end{eqnarray*}
 \end{itemize}
\end{Lem}

Recall the adjoint representation of a Lie $2$-algebra $\g$ on itself as defined in Definition \ref{adjoint}. 
\begin{Def}\label{defi:derivation}
Let $\g$ be a Lie $2$-algebra. A {\bf derivation} of degree $0$ of $\g$ is a $1$-cocycle of $\g$ with respect to the adjoint representation on itself. An {\bf inner derivation} of degree $0$ of $\g$ is a $1$-coboundary of $\g$ with respect to the adjoint representation on itself.
\end{Def}
A notion of a homotopy derivation was introduced in \cite{DL} using the theory of operads. For Lie $2$-algebras, the similarities and differences of these two definitions were clarified in \cite{LLS}.

Denote by $\Der_{0}(\g)$ the set of derivations of degree $0$ of
$\g$.  Define \[\Der_{-1}(\g)=\gl_{-1}(\g):=\Hom(\g_0,\g_{-1}).\] One  obtains a $2$-vector space:
\begin{equation*}
   \Der(\g): \Der_{-1}(\g)\xrightarrow{\bar{d}}\Der_{0}(\g),
\end{equation*}
where $\bar{d}$ is given by
$\bar{d}(\Theta)=(\delta(\Theta),l_{\delta(\Theta)})$, in which
$\delta(\Theta)=(d\circ \Theta,\Theta\circ d)
\in \gl(\g_0)\oplus\gl(\g_{-1})$ and $l_{\delta(\Theta)}\in \Hom(\wedge^2 \g_0,\g_{-1})$ is 

\[l_{\delta(\Theta)}(x,y)={\Theta}[x,y]-[x,\Theta{y}]-[{\Theta}x,y],\qquad \forall x,y\in \g_0.\]
In addition, define
\begin{eqnarray}
\label{eq:bra01}\{(X,l_{X}),\Theta\}&=&[X,\Theta]_C,\\
\label{bra of der}
\{(X,l_X),(Y,l_Y)\}&=&([X,Y]_C,L_X(l_Y)-L_Y(l_X)),
\end{eqnarray}
where $[\cdot,\cdot]_C$ is the commutator bracket and for all $X=(X_0,X_1)\in\gl(\g_0)\oplus\gl(\g_{-1})$, 
$L_X:\Hom(\wedge ^2\g_0,\g_{-1})\longrightarrow\Hom(\wedge^2\g_0,\g_{-1})$  is given by
$$L_X(l_Y)(x,y)=X_1l_Y(x,y)-l_Y(X_0x,y)-l_Y(x,X_0y).$$
\begin{Thm}{\rm(\cite{CSZ})}\label{thm:Der(g)}
With the notations above, $(\Der(\g),\bar{d}, \{\cdot,\cdot\})$ is a strict Lie
$2$-algebra. It is called the {\bf derivation Lie $2$-algebra} of $\g$.
\end{Thm}

Associated with the $2$-term complex $\Der_{-1}(\g)\xrightarrow{\bar{d}} \Der_0(\g)$, there is a 
$3$-term complex of vector spaces
\begin{equation*}
\xymatrix@C=0.5cm{\mathrm{DER}(\mathfrak{g}) : \mathfrak{g}_{-1}\ar[rr]^{d_{\D}}
&& {\Der_{-1}(\mathfrak{g})\oplus{\mathfrak{g}_0}} \ar[rr]^{d_{\D}} &&
               \Der_0(\mathfrak{g}),
              }
\end{equation*}
where the differential $d_\D$ is given by
\begin{eqnarray*}
d_\D(a)&=&([a,\cdot],-da),\qquad \forall a\in \g_{-1};\\
d_\D(\Theta,x)&=&\bar{d}(\Theta)+([x,\cdot],l_3(x,\cdot,\cdot)),\qquad \forall (\Theta,x)\in \Der_{-1}(\g)\oplus \g_0.
\end{eqnarray*}
Moreover, there is also a degree $0$ bracket $[\cdot,\cdot]_\D$ on $\mathrm{DER}(\g)$ given by
\begin{equation*}
\left\{\begin{array}{rcl} [(X,l_X),(Y, l_Y)]_\D&=&\{(X,l_X),(Y,l_Y)\};\\ 
{[(X,l_X),(\Theta,x)]_\D}&=&(\{(X,l_X),\Theta\}+l_X(x,\cdot),X_0x);\\ 
{[(\Theta,x),(\Theta',x')]_\D}&=&-\Theta x'-\Theta' x;\\ 
{[(X,l_X),a]_\D}&=&X_1 a,
\end{array}\right.
\end{equation*}
for any $(X,l_X),(Y,l_Y)\in \Der_0(\g)$, $\Theta,\Theta'\in \Der_{-1}(\g)$, $x,x'\in \g_0$ and $a\in \g_{-1}$, where the bracket $\{\cdot,\cdot\}$ is defined in \eqref{eq:bra01} and \eqref{bra of der}.
\begin{Lem}{\rm(\cite{CSZ})}\label{Der3}
With the above notations, $(\mathrm{DER}(\g), d_\D,
[\cdot,\cdot]_\D)$ is a strict Lie $3$-algebra. It is called the {\bf derivation Lie $3$-algebra} of $\g$.
\end{Lem}

The derivation Lie $2$-algebra and derivation Lie $3$-algebra were introduced in \cite{CSZ} to classify the nonabelian extensions of Lie $2$-algebras.  
\begin{Def}\label{extension}
 \begin{itemize}
 \item[\rm (i)] Let $(\mathfrak{g},d_\g,[\cdot,\cdot]_\g,l_3^\g)$, $(\mathfrak{h},d_\h,[\cdot,\cdot]_\h,l_3^\h)$, $(\hat{\mathfrak{g}},\hat{d},[\cdot,\cdot]_{\hat{\g}},\hat{l_3})$ be Lie $2$-algebras and
$i=(i_{0},i_{1}):\h\longrightarrow\hat{\mathfrak{g}},~~p=(p_{0},p_{1}):\hat{\mathfrak{g}}\longrightarrow\mathfrak{g}$
be strict Lie $2$-algebra homomorphisms. The following sequence of Lie $2$-algebras is a
short exact sequence if $\img i=\ker p$,
$\ker i=0$ and $\img p=\g$:
\begin{equation}\label{eq:ext1}
\CD
  0 @> >>  \h_{-1} @>i_1>> \hat{\g}_{-1} @>p_1>> \g_{-1} @> >> 0 \\
  @. @V d_\h VV @V \hat{d} VV @V d_\g VV @. \\
  0 @> >> \h_{0} @>i_0>> \hat{\g}_0 @>p_0>> \g_0@> >>0.
\endCD
\end{equation}


We call $\hat{\mathfrak{g}}$  an {\bf extension} of $\mathfrak{g}$ by
$\h$, and denote it by $\mathrm{E}_{\hat{\g}}.$
It is called an {\bf abelian extension} if $\h$ is abelian, i.e., $[\cdot,\cdot]_{\h}=0$ and $l^{\h}_3(\cdot,\cdot,\cdot)=0$. 

\item[\rm (ii)] Two extensions of Lie $2$-algebras
 $\mathrm{E}_{\hat{\g}}:0\longrightarrow\h\stackrel{i}{\longrightarrow}\hat{\g}\stackrel{p}{\longrightarrow}\g\longrightarrow0$
 and $\mathrm{E}_{\tilde{\g}}:0\longrightarrow\h\stackrel{j}{\longrightarrow}\tilde{\g}\stackrel{q}{\longrightarrow}\g\longrightarrow0$ are {\bf equivalent},
 if there exists a Lie $2$-algebra homomorphism $F:\hat{\g}\longrightarrow\tilde{\g}$  such that $F\diamond i=j$, $q\diamond
 F=p$ and $F_2(i(u),\alpha)=0$, for all
 $u\in\h_0,~\alpha\in\hat{\g}_0$.
\end{itemize}
\end{Def}

\begin{Thm}\rm{(\cite{CSZ})}\label{nonabelian}
There is a one-to-one correspondence between isomorphism classes of extensions of Lie $2$-algebras given by \eqref{eq:ext1} for $\hat{\g}=\g\oplus \h$ and equivalence classes of morphisms from the Lie $2$-algebra $\g$ to the derivation Lie $3$-algebra $\mathrm{DER}(\h)$.
\end{Thm}

\subsection{The second cohomology and  deformations of Lie $2$-algebras}
The infinitesimal deformation of a Lie $2$-algebra is characterized by the second cohomology of this Lie $2$-algebra with respect to the adjoint representation on itself defined in Example \ref{adjoint}. Moreover, the abelian extension of a Lie $2$-algebra, as a special case of Theorem \ref{nonabelian}, has a simpler description by means of the second cohomology.

Let $(\g, d, [\cdot,\cdot], l_3)$ be a Lie 2-algebra. A {\bf Lie $2$-algebra infinitesimal deformation} of $\g$ is a sequence of linear maps
\[d_i:\g_{-1}\longrightarrow \g_0,\qquad [\cdot,\cdot]_i:\g_0\wedge \g_j\longrightarrow \g_j, j=0,-1,\qquad  l_{3,i}:\wedge^3\g_0\longrightarrow\g_{-1},\qquad i=0,1,\]
with $d_0=d,[\cdot,\cdot]_0=[\cdot,\cdot]$ and $l_{3,0}=l_3$, such that the $\mathbbm{R}[[\lambda]]/(\lambda^2)$-linear operations $(d_\lambda,[\cdot,\cdot]_\lambda, l_3^\lambda)$ on $\g[[\lambda]]/(\lambda^2)$ determined by 
\begin{eqnarray}
\label{inde1} d_\lambda a&:=&d a+\lambda d_1 a,\qquad a\in \g_{-1},\\
\label{inde2} {[x,y]}\dlam&:=& [x,y]+ \lambda [x, y]_1,\\
\label{inde3} {[x,a]}\dlam&:=& [x,a]+ \lambda [x, a]_1,\\
\label{inde4} l_3^\lambda (x, y, z)&:=& l_3(x, y, z)+ \lambda l_{3,1}(x, y, z),\qquad x,y,z\in \g_{0}
 \end{eqnarray}
give a Lie $2$-algebra structure on $\g$ for each $\lambda$. Denote an infinitesimal deformation by $(\g,d_1,[\cdot,\cdot]_1,l_{3,1})$.
Two infinitesimal deformations $\g_\lambda=(\g,d_1,[\cdot,\cdot]_1,l_{3,1})$ and $\g'_\lambda=(\g,d'_1,[\cdot,\cdot]'_1,l'_{3,1})$ of a Lie $2$-algebra $\g$ are said to be {\bf equivalent} if there exists a triple $\phi=(\phi_0,\phi_1,\phi_2)$ with \[\phi_0\in \gl(\g_0),\qquad \phi_1\in \gl(\g_{-1}),\qquad \phi_2\in \Hom(\wedge^2 \g_0,\g_{-1})\] such that 
$\mathrm{Id}+\lambda \phi:\g_\lambda\to \g'_\lambda$ modulo $\lambda^2$ is a   Lie $2$-algebra homomorphism for each $\lambda$.

\begin{Thm}\label{thm:deformation}\rm{(\cite{LSZ})} Let $\g$ be a Lie $2$-algebra. 
There is a one-to-one correspondence between  equivalence classes of Lie $2$-algebra infinitesimal deformations of $\g$ and elements in the second cohomology group $\mathrm{H}^2(\g,\g)$ of $\g$ with respect to the adjoint representation on itself.
\end{Thm}
Let us write down the infinitesimal deformation of $(\g,d,[\cdot,\cdot],l_3)$ given by a $2$-coboundary. Let $(X,l_X)\in C^1(\g,\g)$, with $X=(X_0,X_1)\in \gl(\g_0)\oplus \gl(\g_{-1})$ and $l_X\in \Hom(\wedge^2 \g_0,\g_{-1})$.  Note that 
\[C^2(\g,\g)= \Hom(\g_{-1},\g_0)\oplus \Hom(\wedge^2 \g_0,\g_0)\oplus \Hom(\g_0\wedge \g_{-1},\g_{-1})\oplus \Hom(\wedge^3 \g_0,\g_{-1}).\]
The $2$-coboundary $\D(X,l_X)\in C^2(\g,\g)$  has four components, denoted by  $d_1,[\cdot,\cdot]_1$ (two components) and $l_{3,1}$ respectively. By Lemma \ref{1chain}, we have 
\begin{eqnarray}
\label{inde5}d_1 a&=&dX_1a -X_0 da,\\
\label{inde6} [x,y]_1&=&[x,X_0y]+[X_0x,y]-X_0[x,y]+dl_X(x,y),\\
\label{inde7}[x,a]_1&=&[x,X_1a]+[X_0x,a]-X_1[x,a]+l_X(x,da),\\ 
\label{inde8}l_{3,1}(x,y,z)&=&l_X(x,[y,z])+[x,l_X(y,z)]+l_3(X_0x,y,z)+c.p.-X_1l_3(x,y,z).
\end{eqnarray}
Then the infinitesimal deformation of $\g$ given by the $2$-couboundary $\D(X,l_X)$ is given by formulars \eqref{inde1}-\eqref{inde4} with $d_1,[\cdot,\cdot]_1$ and $l_{3,1}$ given  by the formulas \eqref{inde5}-\eqref{inde8}.

Besides describing the deformations of a Lie $2$-algebra, the second cohomology can also be used to classify abelian extensions of a Lie $2$-algebra, which simplifies and strengthens the result in Theorem \ref{nonabelian} for general extensions.

\begin{Thm}\rm{(\cite{LSZ})}
Given a representation ($\phi_0,\phi_1,\phi_2) : \g\to  \gl(\h)$ of a Lie $2$-algebra $\g$ on a $2$-vector space $\h$, there is a one-to-one correspondence between equivalence classes of abelian extensions of the Lie 2-algebra $\g$ by $\h$ and elements in the second cohomology  group  $\mathrm{H}^2(\g,\h)$.
\end{Thm}

\subsection{The third cohomology and crossed modules of Lie $2$-algebras}
 The third cohomology of a Lie algebra classifies Lie algebra crossed modules due to Gerstenhaber; see Theorem \ref{Gerstenhaber}. This motivated the authors in \cite{LL} to find a notion of crossed modules for Lie $2$-algebras, which is classified by the third cohomology of Lie $2$-algebras.  We shall recall this result in this subsection.

 Let $(\m, d_\m,[\cdot,\cdot]_\m,l^\m_3)$ and $(\g,d_\g,[\cdot,\cdot]_\g,l^\g_3)$ be two Lie 2-algebras.
 We say that {\bf $\g$ acts on $\m$ by derivations} if there exists a Lie $2$-algebra homomorphism
 $$\phi = (\phi_0,\phi_1,\phi_2):\g\longrightarrow \gl(\m)$$
 and a linear map $l_{\phi_0(x)}:\wedge^2\m_0\longrightarrow\m_{-1}$
such that $(\phi_0(x),l_{\phi_0(x)})\in \Der_0(\m)$ and the
map
 $$((\phi_0,l_{\phi_0}),\phi_1,\phi_2):\g\longrightarrow \Der(\m)$$ is a Lie $2$-algebra homomorphism. See Theorem \ref{thm:Der(g)} for the derivation Lie $2$-algebra $\Der(\m)$. 
 Then we shall define a {\bf crossed product} of $\g$ and $\m$ as follows: on $\g\oplus \m$, define $L_1=d_\g+d_\m$ and
\begin{equation*}
\left\{\begin{array}{rcl}
L_2(x+\alpha,y+\beta)&=&[x,y]_\g+[\alpha,\beta]_\m+x\triangleright{\beta}-y\triangleright{\alpha},$\ \ \ \ $\forall{x,y}\in{\mathfrak{g}},\forall{\alpha,\beta}\in{\mathfrak{m}},\\
L_3(x+\alpha,y+\beta,z+\gamma)&=&{l_3^\g(x,y,z)}+l_3^\m(\alpha,\beta,\gamma)
    -(x,y)\triangleright \gamma -(y,z)\triangleright \alpha\\  &&
    -(z,x)\triangleright \beta+l_{\phi_{0}(x)}(\beta,\gamma)+l_{\phi_{0}(y)}(\gamma,\alpha)
    +l_{\phi_{0}(z)}(\alpha,\beta),\\ && \forall{x,y,z}\in{\mathfrak{g}_0},\forall{\alpha,\beta,\gamma}\in{\mathfrak{m}_0}.\end{array}\right.
\end{equation*}
It can be checked that $(\g\oplus \m,L_1,L_2,L_3)$ is  still a Lie $2$-algebra with $\mathfrak{g}$ as a Lie $2$-subalgebra and $\mathfrak{m}$ as an ideal. Denote it by $\g\triangleright_\phi \m$. 
 
 \begin{Def}\rm{(\cite{LL})}
A {\bf crossed module of Lie 2-algebras} is a quadruple
$(\m,\g,\phi,\Pi)$, where $\m$ and $\g$ are two Lie $2$-algebras, $\phi$
is an action of $\g$ on $\m$ by derivations, and
$\Pi:\g\triangleright_\phi \m \rightarrow  \g $ is a Lie $2$-algebra
homomorphism such that $\Pi|_{\g}=\mathrm{Id}_\g=(\mathrm{id}_{\g_0},\mathrm{id}_{\g_{-1}},0)$ and
\begin{itemize}
\item[\rm(i)] $[\alpha,\beta]_\m=\Pi(\alpha)\triangleright{\beta},\ \ \ \ \forall \alpha,\beta\in{\mathfrak{m}}$,
\item[\rm(ii)] $l^\m_3(\alpha,\beta,\gamma)=-(\Pi_0 \alpha,\Pi_0 \beta)\triangleright \gamma
    -\Pi_2(\Pi_0\alpha,\beta)\triangleright{\gamma},\ \ \ \ \forall{\alpha,\beta,\gamma}\in{\mathfrak{m}_0},$
\item[\rm(iii)] $l_{\phi_0(x)}(\beta,\gamma)=-(x,\Pi_0 \beta)\triangleright \gamma
    -\Pi_2(x,\beta)\triangleright{\gamma},\ \ \ \ \forall{\beta,\gamma}\in{\mathfrak{m}_0},x\in{\mathfrak{g}_0}$,
\item[\rm(iv)]
    $\Pi_2(\alpha,\beta)=\Pi_2(\Pi_{0}\alpha,\beta)=\Pi_2(\alpha,\Pi_{0}\beta)
   ,\ \ \ \ \forall{\alpha,\beta}\in{\mathfrak{m}_0}$.
\end{itemize}
In particular, it is called a {\bf strong} crossed module of Lie
$2$-algebras if ${\Pi_2}=0$.
\end{Def}
For a crossed module $(\m,\g,\phi,\Pi)$, decompose $\Pi$
into
$$\Pi=(\Pi_0,\Pi_1,\Pi_2)=\mathrm{Id}_\g + \sigma + \varphi=\big((\mathrm{id}_{\g_0},\varphi_0),(\mathrm{id}_{\g_{-1}},\varphi_1),(0,\sigma,\varphi_2)\big),$$ where $\varphi=\Pi|_{\m}$ and $\sigma={\Pi_2}|_{\g_0\wedge \m_0}$. We see that $\varphi=(\varphi_0,\varphi_1,\varphi_2):\m\longrightarrow \g$ is a Lie 2-algebra homomorphism and $\varphi_2$ is determined by $\sigma$. A crossed module is also denoted by $(\m,\g,\phi,\varphi,\sigma)$. Then a strong crossed module means that $\sigma=0$.



The following example is another motivation for the definition of crossed modules of Lie $2$-algebras.
\begin{Ex}\label{Dcm}  For a Lie $2$-algebra $(\mathfrak{g},d,[\cdot,\cdot],l_{3})$, the derivation Lie $2$-algebra $\Der(\mathfrak{g})$ acts on $\mathfrak{g}$ by derivations with the identity $\mathrm{Id}:\Der(\g)\longrightarrow\Der(\g)$. Consider the extended adjoint representation $\overline{\ad}:\mathfrak{g}\longrightarrow \Der(\mathfrak{g})$ given by 
\begin{equation*}
\overline{\ad}_0(x)=([x,\cdot], l_3(x,\cdot,\cdot)),\qquad\overline{\ad}_1(a)=[a,\cdot],
\qquad \overline{\ad}_2(y,z)=-l_3(y,z,\cdot) 
\end{equation*}
for $x,y,z\in\g_0,a\in \g_{-1}$ and  a linear map $\sigma:\Der_0(\g)\wedge\g_0\longrightarrow \Der_{-1}(\g)$ defined by
\[\sigma((X,l_X),x)=-l_X(x,\cdot), \ \ \ \ \ \forall (X,l_X)\in{\Der_0(\mathfrak{g})}, {x}\in{\mathfrak{g}_{0}}.\]
It can be verified that
$(\g,\Der(\g),\mathrm{Id},\overline{\ad},\sigma)$ is a crossed module of Lie $2$-algebras. 
This crossed module is not strong even if $\g$ is a strict Lie
$2$-algebra. 
\end{Ex}

Associated with a crossed module $(\m,\g,\phi,\Pi)$ of Lie $2$-algebras, there is a strict Lie $3$-algebra on the complex 
\begin{equation*}
\xymatrix@C=0.5cm{ \mathfrak{m}_{-1}\ar[rr]^{d_\D} &&
{\mathfrak{g}_{-1}\oplus{\mathfrak{m}_{0}}} \ar[rr]^{d_\D} &&
                \mathfrak{g}_{0};}
\end{equation*}
see \cite{LL} for details. In particular, for Example \ref{Dcm}, it recovers the derivation Lie $3$-algebra $\mathrm{DER}(\g)$  in Lemma \ref{Der3}.

A crossed module $(\m,\g,\phi,\varphi,\sigma)$ yields a 4-term
exact sequence of 2-vector spaces
\begin{equation*}\label{seq:Lie 2 cm}
\xymatrix@C=0.5cm{ 0 \ar[r] & \mathbb{V }\ar[rr]^{i} &&
                \mathfrak{m} \ar[rr]^{\varphi} && \mathfrak{g} \ar[rr]^{\pi} && \mathfrak{h}\ar[r]  & 0,
                }
\end{equation*}
where $\mathbb{V}:=\ker \varphi, \mathfrak{h}:=\coker
\varphi$, and $i,\pi$ are the canonical inclusion and projection. By
definition, $\mathbb{V}$ is in the
center of $\m$.
If this crossed module is strong, we further get that $\mathfrak{h}$ is a Lie $2$-algebra and there exists an
induced action of $\h$ on $\mathbb{{V}}$.

 Denote by $\mathcal{C}(\mathfrak{h},\mathbb{V})$ the set of strong crossed modules with respect to fixed kernel $\mathbb{V}$, fixed cokernel $\mathfrak{h}$ and fixed action of $\mathfrak{h}$ on $\mathbb{V}$.


\begin{Thm}\rm{(\cite{LL})}
For a Lie $2$-algebra $\mathfrak{h}$ and an $\mathfrak{h}$-module $\mathbb{V}$, there exists a canonical bijection$$\mathcal{C}(\mathfrak{h},\mathbb{V})/_\sim \stackrel{\approx}\longrightarrow\mathrm{H}^{3}(\mathfrak{h},\mathbb{V}).$$
We refer the reader to \cite[Definition 5.8]{LL} for the definition of the equivalence relation $\sim$.\end{Thm}


\section{Integration of Lie $2$-algebras}\label{Integration}

Like Lie $2$-algebras, Lie $2$-groups are the categorification of Lie groups \cite{BL}. In this section, we shall focus on strict Lie $2$-groups and the corresponding strict Lie $2$-algebras. In particular, the automorphisms of a Lie $2$-algebra form a strict Lie $2$-group, whose infinitesimal is the derivation Lie $2$-algebra in Theorem \ref{thm:Der(g)}. 

\subsection{Strict Lie $2$-groups and strict Lie $2$-algebras}\label{Liegpalgcm}

Lie $2$-groups, as the categorification of Lie groups, were introduced and studied in \cite{BL}, where the associative law is replaced by an isomorphism, called the associator, which satisfies  the pentagon equation. We will focus on the strict case.
A strict $2$-group is a group object in Cat. It is a strict $2$-category with one object and  all the $1$-morphisms and $2$-morphisms are invertible. See Appendix A.
Explicitly, a {\bf strict Lie $2$-group} is a Lie groupoid $\Gamma_1\rightrightarrows \Gamma_0$, where $\Gamma_1$ and $\Gamma_0$ are Lie groups and all the structure maps are group homomorphisms. 

We have a Lie group counterpart of the  Lie algebra crossed module defined in Definition \ref{Liecm}. 
\begin{Def}
A {\bf Lie group crossed module} consists of a Lie group homomorphism $\Phi: H\to G$ and an action of $G$ on $H$ by automorphisms (i.e., $\triangleright: G\to \Aut(H)$ is a Lie group homomorphism) satisfying
\begin{eqnarray*}
\Phi(g\triangleright h)&=&g\Phi(h)g^{-1};\\
\Phi(h)\triangleright h'&=&hh'h^{-1},
\end{eqnarray*}
for all $g\in G$ and $h,h'\in H$.
\end{Def}

\begin{Pro}
There is a bijection between strict Lie $2$-groups and crossed modules of Lie groups.
\end{Pro}
Let us sketch the correspondence. If $\Gamma_1\rightrightarrows \Gamma_0$ is a strict Lie $2$-group, taking the restriction of the target map $t$ on $\ker s\subset \Gamma_1$, 
we obtain a Lie group crossed module
 $t: \ker s\to \Gamma_0$, where the action of $\Gamma_0$ on $\ker s$ is $g\triangleright h=\iota(g) h \iota(g)^{-1}$ for $g\in \Gamma_0$ and $h\in \ker s$. Here $\iota:\Gamma_0\to \Gamma_1$ is the inclusion map and the multiplication in $\iota(g) h \iota(g)^{-1}$ is the group multiplication in $\Gamma_1$.

 Conversely, given a Lie group crossed module $\Phi: H\to G$, we build  a strict Lie $2$-group structure on $H\times G\rightrightarrows G$, denoted by $H\rtimes G$, as follows: The group structure on $H\rtimes G$ is given by
 \[(h_1,g_1)*(h_2,g_2)=(h_1(g_1\triangleright h_2),g_1g_2),\qquad \forall g_1,g_2\in G,h_1,h_2\in H,\]
and the groupoid structure on $H\rtimes G$ is:
\begin{itemize}
\item  source and target maps: $s(h,g)=g,t(h,g)=\Phi(h)g$;
\item  inclusion map: $\iota(g)=(e,g)$;
\item  groupoid multiplication $(h_1,g_1)(h_2,g_2)=(h_1h_2,g_2)$ if $g_1=\Phi(h_2)g_2$.
\end{itemize}
The notation $H\rtimes G$  is to emphasize the action of $G$ on $H$ in the crossed module. The Lie group  $H$ also acts on $G$ through $\Phi$, namely, $h\triangleright g=\Phi(h)g$. The Lie groupoid $H\rtimes G\rightrightarrows G$ is actually an action Lie groupoid given by the action of $H$ on $G$.

Various examples of Lie group crossed modules are listed in  Section 8.4 of \cite{BL}. We recall two examples here.


\begin{Ex}
Let $H$ be a Lie group and $\Aut(H)$ the automorphism group. Then we have a Lie group crossed module
\[\Phi:H\to\Aut(H),\qquad \Phi(h)(h')=hh'h^{-1},\qquad \forall h,h'\in H,\]
where $\Aut(H)$ acts on $H$ naturally. The corresponding strict Lie $2$-group is called the {\bf strict automorphism $2$-group} of $H$. Particularly, take $H=\mathrm{SU(2)}$ or the multiplicative group of nonzero quaternions. Then $\Aut(H)=\mathrm{SO(3)}$ and $\Phi:\mathrm{SU(2)}\to \mathrm{SO(3)}$ is the universal covering map. These examples of strict $2$-groups play explicit roles in physics.
\end{Ex}

\begin{Ex}
Suppose that $\Phi:H\to G$ is a surjective homomorphism of Lie groups. Then it is equipped with a Lie group crossed module structure if and only if $\Phi$ is a central extension, meaning that $\ker \Phi$ is contained in the center of $H$. Moreover, when this crossed module exists, it is unique.

Suppose $V$ is a finite dimensional real vector space with a skew-symmetric bilinear form $\omega:V\times V\to \mathbbm{R}$. Let $H=V\oplus \mathbbm{R}$ be the Lie group with the product
\[(v,\alpha)(w,\beta)=(v+w,\alpha+\beta+\omega(v,w)).\]
Then we obtain a Lie group crossed module
\[\Phi:V\oplus \mathbbm{R}\to V,\qquad \Phi(v,\alpha)=v.\]
The corresponding strict $2$-group is called the {\bf Heisenberg 2-group} of $(V,\omega)$.

\end{Ex}

Lie groups have Lie algebras. Taking differentials of the structure maps in a Lie group crossed module, we obtain a Lie algebra crossed module.
\begin{Pro}
There is a unique functor from the category of Lie group crossed modules to the category of Lie algebra crossed modules.
\end{Pro}
So by integrating a strict Lie $2$-algebra, we mean to find a Lie group crossed module whose infinitesimal Lie algebra crossed module is the given strict Lie $2$-algebra.  Actually, there exists a unique $2$-functor from the $2$-category of Lie group crossed modules to that of Lie algebra crossed modules; see \cite[Proposition 45]{BC}. 

The integration of a general Lie $2$-algebra is rather complicated. For the string Lie $2$-algebra in Example \ref{string},
we refer to \cite{BC, BSCS}.  In \cite{BSCS}, the authors constructed an infinite-dimensional Lie $2$-group $\mathcal{P}_1 G$ whose Lie $2$-algebra is equivalent to the string Lie $2$-algebra. The objects of 
$\mathcal{P}_1 G$ are based paths in $G$, while the automorphisms of any object form the level-$1$  Kac-Moody central extension of the loop group $\Omega G$.


More recently, for a Lie $n$-algebra,  the author in \cite{H} provided an explicit construction of its integrating Lie $n$-group. This extends the work in \cite{Ge} for nilpotent $L_\infty$-algebras. In particular, for  Lie algebras, this construction gives the simplicial classifying space of the corresponding simply-connected Lie group. For string Lie $2$-algebras, it yields the model of the simplicial nerve of the string groups \cite{BSCS}.

In the case of strict Lie $2$-algebras, the integration by \cite{H} is proved in \cite{SZ1} to be Morita equivalent to the integration of Lie algebra crossed modules to Lie group crossed modules. As an application, the authors in \cite{SZ1} further integrated a non-strict morphism between Lie algebra crossed modules to a generalized morphism between their corresponding Lie group crossed modules, which includes examples of $2$-term representations up to homotopy, nonabelian extensions and up to homotopy Poisson actions of Lie algebras. An integration of such morphisms was also achieved by the technique of butterflies  \cite{N}.

Another quite interesting case is to integrate semidirect product Lie $2$-algebras in Example \ref{semidirect}. In \cite{SZ3}, the authors integrated such a class of nonstrict Lie $2$-algebras  and they obtained strict Lie $2$-groups in the finite dimensional case by using the butterfly method.

\subsection{Integration of derivations for Lie $2$-algebras}
In this subsection, for a Lie $2$-algebra $(\g,d,[\cdot,\cdot],l_3)$, we give the integration of the derivation Lie $2$-algebra $\Der(\g)$ in Theorem \ref{thm:Der(g)}, which is the automorphism Lie $2$-group $\Aut(\g)$.

To define $\Aut(\g)$, first denote by $\Aut_0(\g)$ the set of Lie $2$-algebra automorphisms of $\g$. See Subsection \ref{2.1} for the definition of automorphisms. It is evident that $\Aut_0(\g)$ with the composition $\diamond$ is a Lie group.
Next, define a multiplication on $\gl_{-1}(\g):=\Hom(\g_0,\g_{-1})$ by
\begin{eqnarray*}\label{eq:multi 1}
\tau\star\tau'=\tau+\tau'+\tau\circ d\circ \tau',\qquad \forall ~\tau,\tau' \in \gl_{-1}(\g).
\end{eqnarray*}
Note that $\star$ satisfies the associative law. So $(\gl_{-1}(\g),\star)$ is a monoid, in which the zero map is the identity element.
$\Aut_{-1}(\g)$ is
defined to be the group of units of $\gl_{-1}(\g)$, which is a Lie group.
Define $\partial:\Aut_{-1}(\g)\longrightarrow \Aut_0(\g)$ by 
\begin{equation*}\label{eq:partial}
  \partial(\tau)=(\mathrm{id}+d\circ \tau,\mathrm{id}+\tau\circ d, l^{\mathrm{Id}}_\tau),\qquad \forall \tau\in\Aut_{-1}(\g),
     \end{equation*}  
     where $l_\tau^{\mathrm{Id}}: \wedge^2 \g_0\to \g_{-1}$ is defined by
  \[l^{\mathrm{Id}}_\tau(x,y)=\tau[x,y]-[x,\tau y]-[\tau x, y]-[\tau x,d \tau y],\qquad \forall x,y\in \g_0.\]
 Based on the fact that $\tau\in \Aut_{-1}(\g)$ if and only if $\mathrm{id}+d\circ \tau\in \mathrm{GL}(\g_0)$, we see that $\partial(\tau)\in \Aut_0(\g)$.
Furthermore, $\Aut_0(\g)$ acts on  $\Aut_{-1}(\g)$ naturally:
\begin{equation*}\label{eq:action}
 A\triangleright \tau=A_1\circ\tau\circ A^{-1}_0\qquad \forall A=(A_0,A_1,A_2)\in\Aut_0(\g),~ \tau\in\Aut_{-1}(\g).
\end{equation*}

\begin{Thm}\rm{(\cite{LLS})}With the notations above,  we have that
\begin{itemize}
\item[\rm{(1)}] the quadruple $(\Aut_{-1}(\g),\Aut_0(\g),\partial,\triangleright)$ is a Lie group crossed module, which is called the {\bf automorphism $2$-group} of $\g$,  and denoted by $\Aut(\g)$.
\item[\rm{(2)}] the differentiation of the Lie group crossed module $(\Aut_{-1}(\g),\Aut_0(\g),\partial,\triangleright)$ is the Lie algebra crossed module $(\Der_{-1}(\g),\Der_0(\g),\bar{d},\phi)$.
\end{itemize}
\end{Thm}

Let us consider the example for the string Lie $2$-algebra $\mathrm{String}(\g): \mathbbm{R}\xrightarrow{0} \g$  defined in Definition \ref{string}, where $\g$ is a semisimple Lie algebra.  
A triple $(X_0,t,l_X)$, with $X_0\in\gl(\g),t\in\mathbbm{R}$ and $l_X\in\wedge^2\g^*$, is a derivation of degree $0$  if and only if
$$
X_0\in\Der(\g),\qquad (\mathfrak{D} l_X)(x,y,z)=tl_3(x,y,z)-l_3(X_0x,y,z)-l_3(x,X_0y,z)-l_3(x,y,X_0z), 
$$
 where $\mathfrak{D}:\wedge^\bullet\g^*\to \wedge^{\bullet+1}\g^*$ is the differential of  $\g$ with coefficients in $\mathbbm{R}$.
 Since $\g$ is semisimple, we have $X_0=\ad_x=[x, \cdot]_\g$ for a unique $x\in \g$. 
By the invariance of the Killing form, we have 
 \[l_3(X_0x,y,z)+l_3(X_0y,z,x)+l_3(X_0z,x,y)=0,\]
 and hence \[(\mathfrak{D} l_X)(x,y,z)=tl_3(x,y,z)=0,\] i.e., $t=0$, as the Cartan 3-form $l_3$ is not exact. 
By the fact that $\mathrm{H}^i(\g)=0, i=1,2$ as  $\g$  is semisimple, there exists a unique $\xi\in \g^*$ such that $l_X=\mathfrak{D}(\xi)$. In summary, we have \[\Der_0(\mathrm{String}(\g))=\{(\ad_x,\mathfrak{D}(\xi)),x\in \g,\xi\in \g^*\}.\]
Moreover, the bracket is
\begin{eqnarray*}
  \{(\ad_x,0,\mathfrak{D}(\xi)), (\ad_y,0,\mathfrak{D}(\eta))\}=(\ad_{[x,y]_\g},0,\mathfrak{D}(\ad^*_x\eta-\ad^*_y\xi)).
\end{eqnarray*}

\begin{Ex} \rm{\cite[Proposition 3.4] {LLS}}\label{pro3.4}
For the string Lie $2$-algebra $\mathrm{String}(\g)$, the derivation Lie $2$-algebra $\Der(\mathrm{String}(\g))$ is as follows:
\begin{itemize}
\item $\Der_0(\mathrm{String}(\g))$ is isomorphic to the semidirect product Lie algebra  $\g\ltimes_{\ad^*}\g^*$;
\item $\Der_{-1}(\mathrm{String}(\g))=\g^* $, which is abelian;
\item  the differential $\bar{d}: \Der_{-1}(\mathrm{String}(\g))\to \Der_0(\mathrm{String}(\g))$ is given by $\bar{d}(\Theta)=(0,-\Theta)$.
\end{itemize}
 \end{Ex}

Following a similar discussion as above,  we have the automorphism $2$-group of the string Lie $2$-algebra.
\begin{Ex}
The automorphism $2$-group $
\Aut(\mathrm{String}(\g))$ of the string Lie $2$-algebra  is 
\begin{itemize}
\item $\Aut_0(\mathrm{String}(\g))$ is isomorphic to the semidirect product $G\ltimes\g^*$, where $G$ is the connected and simply-connected Lie group that integrating $\g$;
\item $\Aut_{-1}(\mathrm{String}(\g))=\g^*$, the abelian Lie group;
\item The map $\partial:
\Aut_{-1}(\mathrm{String}(\g))\rightarrow \Aut_0(\mathrm{String}(\g))$ is given by $\partial(\tau)=(e,-\tau)$, where $e$ is the identity in $G$.
\end{itemize}
\end{Ex}

\appendix
\part*{Appendix}\label{part_appendix}
\addcontentsline{toc}{part}{Appendix}

\section{2-Categories}

Roughly speaking, a  $2$-category consists of objects, $1$-morphisms (also called morphisms) between objects, and $2$-morphisms between morphisms.
The morphisms can be composed along the objects, while the $2$-morphisms can be composed in two different 
directions: along objects-called horizontal composition, and along morphisms-called vertical composition. The composition of morphisms is allowed to be associative only up to coherent associator 2-morphisms. When the associativity holds, it is called a strict $2$-category. See \cite{KS}.

\begin{Def}
A {\bf strict $2$-category} $\mathcal{C}$ consists of the following data
\begin{itemize}
\item[\rm (1)] a set of {\bf objects} $\mathrm{Ob}(\C)$;
\item[\rm (2)] for each pair $x,y\in \mathrm{Ob}(\C)$, a category $\Mor_{\C}(x,y)$; 
\begin{itemize}
\item The objects of $\Mor_{\C}(x,y)$ are called 
{\bf $1$-morphisms}  and denoted by $F:x\to y$. 
\item the morphisms between these $1$-morphisms are called 
{\bf $2$-morphisms} and denoted by $\theta:F'\Rightarrow F$. 
\item The composition of $2$-morphisms in $\Mor_{\C}(x,y)$ is called {\bf vertical composition} and denoted as $\theta'\circ \theta$ for $\theta':F'\Rightarrow F$ and $\theta:F''\Rightarrow F'$:
\begin{eqnarray*}
\xymatrix@C+2em{
   \bullet & \bullet
    \ar@/_1.5pc/[l]_-{F''}^{}="0"
    \ar[l]_{F'\qquad}_{}="1"
    \ar@/^1.5pc/[l]^-{F}_{}="2"
    \ar@{=>} "0";"1"^{\theta}
    \ar@{=>} "1";"2"^{\theta'}
  }&=&\xymatrix@C+2em{
  \bullet &
  \ar@/_1pc/[l]_{F''}_{}="0"
  \ar@/^1pc/[l]^{F}="1"
  \ar@{=>}"0";"1"^{\theta'\circ \theta}
  \bullet};\end{eqnarray*}
\end{itemize}
\item[\rm (3)] For each triple $x,y,z\in \mathrm{Ob}(\C)$, a functor 
\[(\circ, *):\Mor_{\C}(y,z)\times \Mor_{\C}(x,y)\to \Mor_{\C}(x,z),\]
the image of the pair of $1$-morphisms $(F,G)$ is called the {\bf composition} of $F$ and $G$, denoted by $F\circ G$. The image of the pair of $2$-morphisms $(\theta,\tau)$ is called the {\bf horizontal composition} and denoted by $\theta*\tau$:
\begin{eqnarray*}
\xymatrix@C+2em{
 \bullet &
  \ar@/_1pc/[l]_{F}_{}="0"
  \ar@/^1pc/[l]^{F'}^{}="1"
  \ar@{=>}"0";"1"^{\theta}
  \bullet &
  \ar@/_1pc/[l]_{G}_{}="2"
  \ar@/^1pc/[l]^{G'}^{}="3"
  \ar@{=>}"2";"3"^{\tau}
  \bullet}
=\xymatrix@C+2em{
 \bullet &
  \ar@/_1pc/[l]_{F\circ G}_{}="0"
  \ar@/^1pc/[l]^{F'\circ G' }="1"
  \ar@{=>}"0";"1"^{\theta*\tau}
  \bullet }.
\end{eqnarray*}
\end{itemize}
These data satisfy the following rules:
\begin{itemize}
\item[\rm (1)] The sets of objects and  $1$-morphisms  with composition of $1$-morphisms form a category;
\item[\rm (2)] The horizontal composition of $2$-morphisms is associative;
\item[\rm (3)] The identity $2$-morphism $\mathrm{id}_{\mathrm{id}_x}$ of the identity $1$-morphism $\mathrm{id}_x$ is a unit for the horizontal composition.
\end{itemize}
\end{Def}

Since $(\circ, *):\Mor_{\C}(y,z)\times \Mor_{\C}(x,y)\to \Mor_{\C}(x,z)$ is a functor, we have the interchange law:
\[(\theta'\circ \theta)*(\tau'\circ \tau)=(\theta'* \tau')\circ (\theta*\tau),\]
for $\theta:F\Rightarrow F', \theta':F'\Rightarrow F'',\tau:G\Rightarrow G'$ and $\tau':G'\Rightarrow G''$.

A strict Lie $2$-group is a strict $2$-category with one object and all the $1$-morphisms and $2$-morphisms are invertible.

\begin{Def}\label{2-homo} Let $\g$ and $\h$ be two $2$-term $L_\infty$-algebras and $\phi,\psi: \g\to \h$ be two $2$-term $L_\infty$-algebra homomorphisms. A {\bf  $2$-homomorphism}  $\tau:\phi\Rightarrow \psi$  is a linear map $\tau: \g_0\to \h_{-1}$ such that 
\begin{itemize}
\item[\rm{(1)}] $\tau$ is a homotopy between the chain maps $(\phi_0,\phi_1)$ and $(\psi_0,\psi_1)$, i.e.,
\[\psi_0-\phi_0=d_{\h}\circ \tau,\qquad \psi_1-\phi_1=\tau\circ d_\g.\]
\item[\rm{(2)}] $\psi_2(x,y)-\phi_2(x,y)=\tau([x,y]_\g)-[\phi_0(x),\tau(y)]_\h-[\tau(x),\psi_0(y)]_\h$ for $x,y\in \g_0$.
\end{itemize}
\end{Def}

\begin{Def}
A {\bf homotopy equivalence} between two $2$-term $L_\infty$-algebras $\g$ and $\h$ is a pair of $L_\infty$-algebra homomorphisms $\phi:\g\to \h$ and $\psi:\h\to \g$ such that the compositions $\phi\diamond \psi$ and $\psi\diamond \phi$ are homotopic to the identity map, i.e., there are $2$-homomorphisms  $\theta: \phi\diamond \psi\Rightarrow \mathrm{Id}_\h$ and $\tau: \psi\diamond \phi\Rightarrow \mathrm{Id}_\g$. 
\end{Def}
Homotopy equivalences for strict graded Lie $2$-algebras we used in Theorem \ref{BCLXth} are similar to be defined. We refer to the Appendix $A$ in \cite{BCLX} for details.

\begin{Ex}
There is a strict $2$-category $\mathrm{2TermL_\infty}$ with $2$-term $L_\infty$-algebras as objects, homomorphisms as $1$-morphisms and $2$-homomorphisms as $2$-morphisms. See Definitions \ref{homo} and \ref{2-homo} for homomorphisms and $2$-homomorphisms respectively. 
\end{Ex}

Generally, let $K$ be a category. There is a strict $2$-category, denoted by $\mathrm{KCat}$,  with categories in $K$ as objects, functors in $K$ as $1$-morphisms and natural transformations in $K$ as $2$-morphisms; see \cite[Proposition 4]{BC}. For instance, when $K$ is the category of vector spaces, $\mathrm{KCat}$ is the strict $2$-category $\mathrm{2Vect}$ of $2$-vector spaces we mentioned in Subsection \ref{2.1}.

\section{Lie algebroids and Lie grouppods}

We refer to Mackenzie's book \cite{Mackenzie} for a thorough description of the theories of Lie algebroids and groupoids.
\begin{Def}
  A {\bf Lie algebroid} is a vector bundle $A\rightarrow M$  with a Lie algebra structure $[\cdot,\cdot]_A$ on
the section space $\Gamma(A)$ and a vector bundle morphism
$\rho_A:A\rightarrow TM$ from $A$ to the tangent bundle $TM$, called the anchor, such that 
$$[u,fv]_A=f[u,v]_A+\rho_A(u)(f)v,\quad \forall u,v\in\Gamma(A),~f\in
C^\infty(M).$$
\end{Def}

A groupoid is a small category such that every arrow is invertible. Explicitly,
\begin{Def}
A {\bf groupoid} is a pair $(\mathcal{G},M)$, where $M$ is the set of objects and $\G$ is the set of arrows, with the  structure maps
\begin{itemize}
\item two surjection maps $s,t: \G\to M$, called the source map and target map, respectively;
\item  the multiplication $\cdot:\G^{(2)}\to \G$, where $\G^{(2)}=\{(g_1,g_2)\in \G\times \G| s(g_1)=t(g_2)\}$;
\item  the inverse map $(\cdot)^{-1}:\G\to \G$;
\item the inclusion map $\iota: M\to \G$, called the identity map,
\end{itemize}
satisfying the following properties:
\begin{itemize}
\item[\rm{(1)}] \rm{(associativity)} $(g_1\cdot g_2)\cdot g_3=g_1\cdot (g_2\cdot g_3)$, whenever the multiplications are well-defined;
\item[\rm{(2)}] \rm{(unitality)} $\iota(t(g))\cdot g=g=g\cdot \iota(s(g))$;
\item[\rm{(3)}]  \rm{(invertibility)} $g\cdot g^{-1}=\iota(t(g))$, $g^{-1}\cdot g=\iota(s(g))$.
\end{itemize}
Denote a groupoid by $(\G\rightrightarrows M,s,t)$ or simply by $\G$.

A {\bf Lie groupoid} is a groupoid such that both the set of objects and the set of arrows  are smooth manifolds, all structure maps are smooth, and source and target maps are subjective submersions.
\end{Def}
As the tangent space of a Lie group at the identity has a Lie algebra structure, for a Lie groupoid $(\mathcal{G}\rightrightarrows M,s,t)$, on the vector bundle $A:=\ker(s_*)|_M\rightarrow M$, there is a Lie algebroid structure defined as follows:
the anchor map $A\rightarrow TM$ is simply $t_*$ and the Lie bracket $[u,v]_A$ is determined by
\[\overrightarrow{[u,v]_A}=[\overrightarrow{u},\overrightarrow{v}],\qquad \forall u,v\in \Gamma(A),\]
where $\overrightarrow{u}$ denotes the right-invariant vector field on $\G$ given by $\overrightarrow{u}_g=R_{g_*}u_{t(g)}$
and $[\cdot,\cdot]$ is the Schouten bracket on $\mathfrak{X}(\mathcal{G})$.

Unlike the Lie algebra case, not every Lie algebroid can be integrated to a Lie groupoid. We refer to \cite{CFAnal} for the integrability condition.
\begin{Ex}
A Lie algebra is a Lie algebroid with the base manifold being a point; A Lie group is a Lie groupoid with the base manifold being a point.
\end{Ex}
\begin{Ex}
For a manifold $M$, the tangent bundle $TM$ with the Schouten bracket is a Lie algebroid whose anchor is the identity map. The Lie groupoid of this Lie algebroid is the pair groupoid $M\times M\rightrightarrows M$, where the source and target maps are the two projections, the inclusion map is $\iota(x)=(x,x)$, and the multiplication is \[(x,y)\cdot (y,z)=(x,z),\qquad \forall  x,y,z\in M.\]

Another Lie groupoid integrating $TM$ is the fundamental groupoid, or homotopy groupoid: 
\[\pi_1(M):=\{paths~in~M\}/homotopies\rightrightarrows M, \]whose source and target are end points of a path, and multiplication is the concatenation of paths. This is the unique source simply-connected groupoid that integrating $TM$; see \cite{CFAnal, Mackenzie} for details.
\end{Ex}
\begin{Ex}
For a Poisson manifold $(M,\pi)$, there is a Lie algebroid structure on the cotangent bundle $T^*M$, with the anchor \[\pi^\sharp: T^*M\to TM,\qquad \pi^\sharp(\xi)(\eta):=\pi(\xi,\eta),\qquad \forall \xi,\eta\in \Omega^1(M),\]
and the Lie bracket
\[[\xi,\eta]_{T^*M}=L_{\pi^\sharp(\xi)} \eta-L_{\pi^\sharp(\eta)} \xi-d_{dR} \pi(\xi,\eta).\]
The Lie groupoid of this Lie algebroid is the symplectic groupoid of the Poisson manifold $(M,\pi)$.
\end{Ex}
\begin{Ex}\label{actiongroupoid}
Let  $\g$ be a Lie algebra which acts on a manifold $M$ by the Lie algebra homomorphism $\phi:\g\to \mathfrak{X}(M)$. Then there is a Lie algebroid structure on the trivial bundle $A:=\g\times M\to M$, where the anchor $\rho_A$ 
is 
\[\rho_A(fu)=f\phi(u),\]
and the Lie bracket $[\cdot,\cdot]_A$ is given by
\[ [fu,gv]_A=fg[u,v]_\g+f\phi(u)(g)v-g\phi(v)(f)u,\qquad \forall f,g\in C^\infty(M),u,v\in \g.\]
This Lie algebroid is called the {\bf action Lie algebroid}.
The corresponding Lie groupoid is the action Lie groupoid given as follows. Let  $G$ be a Lie group which acts a manifold $M$ by $\Phi: G\to \Diff(M)$. Then we have a Lie groupoid structure on $G\times M\rightrightarrows M$, whose source and target maps are 
\[s(g,m)=m,\qquad t(g,m)=\Phi(g)(m),\]
the inclusion map is $\iota(m)=(e,m)$, and the multiplication is 
\[(g,m)\cdot (h,n)=(gh,n),\qquad \forall m,n\in M,g,h\in G,\]
whenever $m=\Phi(h)(n)$.

\end{Ex}

\end{document}